\documentclass[12pt,letterpaper]{article}
\usepackage{amsfonts}
\usepackage{amssymb}
\usepackage{latexsym}  
\usepackage{pifont}  

\usepackage{amsmath}   

\DeclareMathAlphabet{\mathbbold}{U}{bbold}{m}{n}

\let\SavedRightarrow=\Rightarrow
\usepackage{marvosym}
\let\Rightarrow=\SavedRightarrow

\newcommand{\nameX}{\raisebox{0pt}[12pt]{$\overset{\kern1pt\raisebox{-3pt}{\scriptsize $\diamond$}}{X}$}}
\newcommand{\nameF}{\raisebox{0pt}[12pt]{$\overset{\kern1pt\raisebox{-3pt}{\scriptsize $\diamond$}}{F}$}}

\newcommand\KK{{\mathcal K}}
\newcommand\FF{{\mathcal F}}
\newcommand\EE{{\mathcal E}}
\newcommand\BB{{\mathcal B}}
\newcommand\HH{{\mathcal H}}

\newcommand\PP{{\mathcal P}}
\newcommand\II{{\mathcal I}}

\newcommand\RRR{{\mathbb R}}

\newcommand\TTT{{\mathbb T}}
\newcommand\QQQ{{\mathbb Q}}
\newcommand\NNN{{\mathbb N}}
\newcommand\EEE{{\mathbb E}}

\newcommand\cccc{{\mathfrak c}}
\newcommand\KKKK{{\mathfrak K}}

\newcommand\cchi{{\raise 2 pt \hbox{$\chi$}}}

\newcommand{\MS}{\mathsf{MS}}  
\newcommand{\ma}{\mathsf{ma}}  
\newcommand\ro{\mathsf{ro}}  

\newcommand\mosim{\mathord{\sim}}  

\newcommand\one{\mathbbold{1}} 
\newcommand\zero{\mathbbold{0}} 

\newcommand\res{\mathord {\upharpoonright}}  

\newcommand\onto{\twoheadrightarrow}  

\newcommand\MA{\mathrm{MA}}
\newcommand\CH{\mathrm{CH}}
\newcommand\Ldim{\mathrm{Ldim}}   

\newcommand\seq[2]{[#1]^{#2\raisebox{1pt}{\mbox{\tiny$\uparrow$}}}} 

\newcommand\ran{\mathrm{ran}}  
\newcommand\dom{\mathrm{dom}}  
\newcommand\diam{\mathrm{diam}}   
\newcommand\cl{\mathrm{cl}}   

\newcommand\lh{\mathrm{lh}}   

\newcommand\down{{\fam0 \downarrow}}

\newcommand\cat{^{\mathord{\frown}}}  

\newcommand\iv{^{-1}} 

\def\eop{{\Large \Coffeecup}}  

{\end{enumerate}}
\newenvironment{itemizz}{\begin{itemize}\setlength{\itemsep}{-1mm}} %
{\end{itemize}}                              

\newenvironment{itemizn}[1] 
{\begin{itemize} \setlength{\itemsep}{-1mm} %
\renewcommand\labelitemi{\ding{#1}}} %
{\end{itemize}}

\newtheorem{theorem}{Theorem}[section]
\newtheorem{definition}[theorem]{Definition}
\newtheorem{lemma}[theorem]{Lemma}
\newtheorem{corollary}[theorem]{Corollary}
\newtheorem{proposition}[theorem]{Proposition}
\newtheorem{example}[theorem]{Example}

\newenvironment{proof}{{\bf Proof.}}{\eop\medskip}
\newenvironment{proofof}[1]{\medskip \textbf{Proof of #1.}}{\eop\medskip}

\begin{document}

\title{Dissipated Compacta\footnote{
2000 Mathematics Subject Classification:
Primary  54D30, 54F05.  
Key Words and Phrases: scattered space, LOTS, inverse limit,
irreducible map.
}}

\author{
Kenneth Kunen\footnote{University of Wisconsin,  Madison, WI  53706, U.S.A.,
\ \ kunen@math.wisc.edu}
\thanks{Author partially supported by NSF Grant DMS-0456653.}
}

\maketitle

\begin{abstract}
The dissipated spaces form a class of compacta which
contains both the scattered compacta and the compact LOTSes
(linearly ordered topological spaces),
and a number of theorems true for these latter two classes
are true more generally for the dissipated spaces.
For example, every regular Borel measure on a dissipated space
is separable.

A product of two compact LOTSes is usually not dissipated, but it may
satisfy a weakening of that property.  In fact, the degree
of dissipation of a space can be used to distinguish topologically
a product of $n$ LOTSes from a product of $m$ LOTSes.

\end{abstract}

\section{Introduction} 
\label{sec-intro}
All topologies discussed in this paper are assumed to be Hausdorff.
As usual, a subset of a space is \emph{perfect} iff it is closed and
non-empty and has no isolated points, so $X$ is \emph{scattered} iff
$X$ has no perfect subsets.

There are many constructions in the literature which build a
compactum $X$ as an inverse limit of metric compacta $X_\alpha$ for
$\alpha < \omega_1$, with
the bonding maps $\pi^\beta_\alpha : X_\beta \onto X_\alpha$
for $\alpha < \beta < \omega_1$.
In some cases, as in \cite{Fed1, HK2, HK3}, the construction
has the property that for each $\alpha,\beta$,
$(\pi^\beta_\alpha)\iv\{x\}$ is a singleton for all but
countably many $x \in X_\alpha$.
We shall call such $\pi^\beta_\alpha$ 
\emph{tight maps}; these are discussed in greater detail in
Section \ref{sec-tight}.  The spaces $X$ so constructed
are examples of \emph{dissipated compacta}; these are discussed
in Section \ref{sec-dis}.
Section \ref{sec-ab} shows that the property of tightness is
absolute for transitive models of set theory.

The precise definition of ``dissipated'' in  Section \ref{sec-dis} will be
that there are ``sufficiently many'' tight maps onto metric compacta;
so the definition will not mention inverse limits.
Then, Section \ref{sec-inv} will relate this definition to inverse limits.

Dissipated compacta include the scattered compacta,
the metric compacta, and the compact LOTSes
(totally ordered spaces with the order topology).
Section \ref{sec-dis} also describes 
the more general notion of $\kappa$--dissipated, which 
gets weaker as $\kappa$ gets bigger; ``dissipated'' is the same as
as ``$2$--dissipated'', while 
``$1$--dissipated'' is the same as ``scattered''.
Every regular Borel measure on a
$2^{\aleph_0}$--dissipated compactum is separable
(see Section \ref{sec-meas}).

If $X$ is the double arrow space of Alexandroff and  Urysohn,
then $X$ is a non-scattered LOTS and hence is
$2$--dissipated but not $1$--dissipated,
while $X^{n+1}$ is 
$(2^n + 1)$--dissipated but not $2^n$--dissipated.
Considerations of this sort can be used to
distinguish topologically
a product of $n$ LOTSes from a product of $m$ LOTSes;
see Section  \ref{sec-LOTS}.

\section{Tight Maps}
\label{sec-tight}
As usual,
$f : X \to Y$ means that $f$ is a \emph{continuous} map from $X$ to $Y$,
and $f : X \onto Y$ means that $f$ is a continuous map from $X$
\emph{onto} $Y$.

\begin{definition}
\label{def-tight}
Assume that $X,Y$ are compact and $f: X \to Y$.
\begin{itemizn}{43}
\item A \emph{loose family} for $f$ is a disjoint family $\PP$ of closed
subsets of $X$ such that for some non-scattered $Q \subseteq Y$,
$Q = f(P)$ for all $P \in \PP$.
\item $f$ is $\kappa$--\emph{tight} iff there are no loose families 
for $f$ of size $\kappa$.
\item $f$ is \emph{tight} iff $f$ is $2$--tight.
\end{itemizn}
\end{definition}

This notion gets weaker as $\kappa$ gets bigger.
$f$ is $1$--tight iff $f(X)$ is scattered,
so that ``$2$--tight'' is the first non-trivial case.
$f$ is trivially $|X|^+$--tight. 
The usual projection from $[0,1]^2$ onto $[0,1]$ is not
$2^{\aleph_0}$--tight.

Some easy equivalents to ``$\kappa$--tight'':

\begin{lemma}
\label{lemma-tight-equiv}
Assume that $X,Y$ are compact and $f : X \to Y$.  Then
$(1) \leftrightarrow (2)$.
If $\kappa$ is finite, then  $(1) \leftrightarrow (3)$;
if also $Y$ is metric, then all five of the following are equivalent:

\begin{itemizz}
\item[1.] There is a loose family of size $\kappa$.
\item[2.]
There is a disjoint family $\PP$ of perfect subsets of $X$ with $|\PP| = \kappa$
and a perfect $Q \subseteq Y$ such that $Q = f(P)$ for all $P \in \PP$.
\item[3.] There are distinct $a_i\in X$ for $i < \kappa$ with
all $f(a_i) = b \in Y$ such that whenever $U_i$ is a neighborhood
of $a_i$ for $i < \kappa$,
$\bigcap_{i < \kappa} f(\overline{U_i})$ is not scattered.
\item[4.] For some metric $M$ and $\varphi\in C(X,M)$,
$\{y \in Y : |\varphi(f\iv\{y\})| \ge \kappa\}$ is uncountable.
\item[5.] Statement $(4)$, with $M = [0,1]$.
\end{itemizz}
\end{lemma}
\begin{proof}
$(2) \to (1)$ is obvious.
Now, assume (1), and let $\PP$ be a loose family of size $\kappa$,
with $Q = f(P)$ for $P \in \PP$.  Let $Q'$ be a perfect subset of $Q$,
and, for $P \in \PP$, let $P'$ be a closed subset of
$P \cap f\iv(Q')$ such that $f \res P' : P' \onto Q'$ is irreducible.
Then $\{P': P \in \PP\}$ satisfies (2).

From now on assume that $\kappa$ is finite.

$(3) \to (1)$ and $(5) \to (4)$ are obvious.

For $(1) \to (3)$, use compactness of $\prod_i P_i$ 
and the fact that a finite union of scattered spaces is scattered.

For $(1) \to (5)$:  
If $\PP = \{P_i : i < \kappa\}$ is a loose family, with $Q = f(P_i)$,
apply the Tietze Theorem to get
$\varphi\in C(X, [0,1])$ such that $\varphi(x) = 2^{-i}$ for all
$x \in P_i$.

Now, we prove $(4) \to (1)$ when $Y$ is metric.
Fix $\varphi$ as in $(4)$.  We may assume that $M = \varphi(X)$,
so that $M$ is compact.  Let $\BB$ be a countable base for $M$.
Then we can find $B_i \in \BB$ for $i < \kappa$ such that the
$\overline{B_i}$ are disjoint and such that $Q :=
\{y \in \nobreak Y :\forall i < \kappa \, [\varphi(f\iv\{y\}) \cap
\overline{B_i} \ne \emptyset]\}$
is uncountable, and hence not scattered (since $Y$ is metric).
$Q$ is also closed.  Let $P_i =  f\iv(Q) \cap \varphi\iv(\overline{B_i})$.
Then $\{P_i : i < \kappa\}$ is a loose family.
\end{proof}

\begin{lemma}
\label{lemma-tight-LOTS}
If $X,Y$ are compact LOTSes and $f:X \to Y$ is order-preserving
$(x_1 < x_2 \to f(x_1) \le f(x_2))$,
then $f$ is tight.
\end{lemma}
\begin{proof}
If not, we would have $a_0 < a_1$ and $b$
as in (3) of Lemma \ref{lemma-tight-equiv}.
Let $U_0, U_1$ be open intervals in $X$ with disjoint closures
such that each $a_i \in U_i$.  But then
$f(\overline{U_0}) \cap f(\overline{U_1})  = \{b\}$, a contradiction.
\end{proof}

In many cases, the loose family will be defined uniformly
via a continuous function, and we may replace the cardinal $\kappa$
in Definition \ref{def-tight} by some compact space $K$ of size $\kappa$:

\begin{definition}
\label{def-K-tight}
Assume that $X,Y,K$ are compact spaces and $f: X \to Y$.
Then a \emph{$K$--loose function} for
$f$ is a $\varphi : \dom(\varphi) \to K$ such that\textup: $\dom(\varphi)$ is 
closed in $X$, and for some non-scattered $Q \subseteq Y$,
$\varphi(f\iv\{y\}) = K$ for all $y \in Q$.
\end{definition}

Note that we then have a loose family $\PP = \{P_z: z \in K\}$
of size $|K|$, where $P_z = f\iv(Q) \cap \varphi\iv\{z\} $.
For finite $n$, we may view the ordinal $n$ as a discrete topological
space, so an $n$--loose \textit{function} is equivalent to a loose
\textit{family} $\PP = \{P_i : i < n\}$, since $\varphi$ can map
$P_i$ to $i \in n$.  The same phenomenon holds for $\aleph_0$, but
seems harder to prove:

\begin{theorem}
\label{thm-omega-loose}
If $X,Y$ are compact and $f: X \to Y$, then there is an infinite loose family
iff there is an $(\omega + 1)$--loose function.
\end{theorem}

This will be proved in Section \ref{sec-ab}.
Beyond $\aleph_0$, there is no simple equivalence between the
cardinal version and the topological version of looseness.
At $2^{\aleph_0}$, we shall use the following terminology
to avoid possible confusion between
the Cantor set $2^\omega$ and the cardinal $\cccc = 2^{\aleph_0}$:

\begin{definition}
\label{def-weak-tight}
Assume that $X,Y$ are compact and $f: X \to Y$.
\begin{itemizn}{43}
\item A \emph{strongly $\cccc$--loose family} for
$f$ is a $K$--loose function  $\varphi : \dom(\varphi) \to K$,
where $K$ is the Cantor set $2^\omega$.
\item $f$ is \emph{weakly $\cccc$--tight} iff there is no
strongly $\cccc$--loose function for $f$.
\end{itemizn}
\end{definition}

In this paper, whenever we produce a loose family
of size $2^{\aleph_0}$, it will usually be strongly $\cccc$--loose.
However, if we view $\cccc + 1$ as a compact ordinal and let
$X = Y \times (\cccc + 1)$, then assuming that $Y$ is not scattered,
the usual projection $f: X \onto Y$ has an obvious
loose family of size $\cccc$ but no strongly $\cccc$--loose family.

When $X,Y$ are both metric, the $\kappa$-tightness of $f$ is related to
the sizes of the sets $f\iv\{y\}$ by:

\begin{theorem}
\label{thm-tight-metric}
If $X,Y$ are compact metric and $f: X \to Y$,
then $f$ is $\kappa$--tight iff
$\{y \in Y : |f\iv\{y\}| \ge \kappa\}$ is countable.
$f$ is weakly $\cccc$--tight iff $f$ is $\cccc$--tight.
\end{theorem}

In particular, if $f : X \onto Y$, then $f$ is tight iff
$f\iv\{y\}$ is a singleton for all but
countably many $y$, as we said in the Introduction.

For both ``iff''s, the $\leftarrow$ direction is trivial and is true for
any $X,Y$. 
For $\kappa = 3$, say, the proof of the 
$\rightarrow$ direction will show that if there are
uncountably many $y \in Y$ such that $f\iv\{y\}$ contains three or more
points, then for some perfect $Q \subseteq Y$, we can, on $Q$,
choose three of these points continuously, producing
disjoint perfect $P_0,P_1,P_2 \subseteq X$ which
$f$ maps homeomorphically onto $Q$,
so $\{ P_0,P_1,P_2\}$ is a loose family of size $3$.

Since $X$ is second countable, each $f\iv\{y\}$ is either countable
or of size $2^{\aleph_0}$, so it is sufficient to prove the theorem
for the cases $\kappa \le \aleph_0$ and $\kappa = 2^{\aleph_0}$.
However, for $\kappa = 2^{\aleph_0}$, we can get more detailed results.
For example, if there are
uncountably many $y \in Y$ such that $f\iv\{y\}$ contains a Klein bottle,
then we can choose the bottle continuously on a perfect set
(see Theorem \ref{thm-metric-embed}).
This ``continuous selector'' result follows easily from standard
descriptive set theory.  First, observe:

\begin{lemma}
\label{lemma-suslin}
Suppose that $g : \Phi \to Y$, where $Y$ is a Polish space,
$\Phi$ is an analytic subset of some Polish space, and $g(\Phi)$
is uncountable.  Then there is a Cantor subset $C \subseteq \Phi$
such that $g$ is 1-1 on $C$.
\end{lemma}
\begin{proof}
Let $h : \omega^\omega \onto \Phi$, apply the classical argument
of Suslin to obtain a Cantor subset $D \subseteq \omega^\omega $ such that 
$g\circ h$ is 1-1 on $D$, and let $C = h(D)$.
\end{proof}

\begin{theorem}
\label{thm-metric-embed}
Assume that $X,Y,Z$ are compact metric, $f: X \to Y$,
and there are uncountably many $y \in Y$ such that $f\iv\{y\}$
contains a homeomorphic copy of $Z$. 
Then there is a perfect $Q \subseteq Y$ and a 1-1 map $i : Q \times Z \to X$
such that $f(i(q,u)) = q$ for all $(q,u) \in Q\times Z$.
\end{theorem}
\begin{proof}
Assume that $Z \ne \emptyset$.
Fix metrics $d_Z, d_X$ on $Z,X$, and give $C(Z,X)$ the usual
uniform metric, which makes it a Polish space.
Let $\Phi$ be the set of all $\varphi \in C(Z,X)$ such that
$\varphi$ is 1-1 and $\varphi(Z) \subseteq f\iv\{y\}$ for some
(unique) $y \in Y$.
Observe that $\Phi$ is an $F_{\sigma\delta}$ set, since the
``$\varphi$ is 1-1'' can be expressed as:
\[
\forall \varepsilon > 0 \, \exists \delta > 0 \, \forall u,v \in Z\,
[d_Z(u,v) \ge \varepsilon \to d_X(\varphi(u), \varphi(v)) \ge \delta] \ \ .
\]
Define $g : \Phi \to Y$ so that $g(\varphi)$ is the $y \in Y$
such that $\varphi(Z) \subseteq f\iv\{y\}$.
Using Lemma \ref{lemma-suslin}, let $C \subseteq \Phi$ be a
Cantor subset with $g$ 1-1 on $C$, let $Q = g(C)$,
and let $i( g(\varphi), u) = \varphi(u)$.
\end{proof}

\begin{proofof}{Theorem \ref{thm-tight-metric}}
To prove the $\rightarrow$ direction of the first ``iff'' in the three cases 
$\kappa < \aleph_0$, $\kappa = \aleph_0$, and $\kappa = \cccc$, 
apply Theorem \ref{thm-metric-embed} respectively with $Z$ the space
$\kappa$ (with the discrete topology), 
$\omega+1$, and $2^\omega$.
This also yields the $\rightarrow$ direction of the second ``iff''.
\end{proofof}

Of course, we are using the fact that every uncountable
metric compactum contains a copy of the Cantor set.
One could also prove Theorem \ref{thm-tight-metric} using the following,
plus the fact that every uncountable metric compactum maps onto $[0,1]$:

\begin{theorem}
\label{thm-metric-cover}
Assume that $X,Y,K$ are compact metric with $f: X \to Y$,
and assume that for uncountably many $y \in Y$,
there is a closed subset of $f\iv\{y\}$ which can be mapped onto $K$.
Then there is a $K$--loose function for $f$.
\end{theorem}
\begin{proof}
Let $H$ be the Hilbert cube, $[0,1]^\omega$. 
We may assume that $K \subseteq H$.  Then, for uncountably many $y \in Y$,
there is a $\psi \in C(X,H)$ such that $\psi(f\iv\{y\}) \supseteq K$.
Let $\Psi = \{(y, \psi) \in Y \times C(X,H) : \psi(f\iv\{y\}) \supseteq K \}$,
and let $g(y, \psi) = y$.  Applying Lemma \ref{lemma-suslin},
let $C \subseteq \Psi$ be a Cantor set on which $g$ is 1-1,
and let $Q = g(C) \subseteq Y$.
For $(y,\psi) \in C$, let $E_y = \{x \in X : \psi(x) \in K\}$.
Define $\varphi$ so that $\dom(\varphi) = \bigcup\{E_y : y \in Q\}$,
and $\varphi(x) = \psi(x)$ whenever $x \in \dom(\varphi)$
and $(y,\psi) \in C$.  Then $\varphi$ is a $K$--loose function.
\end{proof}

Theorems  \ref{thm-tight-metric}, \ref{thm-metric-embed},
and \ref{thm-metric-cover}
can fail when $X$ is not metric; counter-examples are provided
by the double arrow space and some related spaces described by:

\begin{definition}
\label{def-das}
$I = [0,1]$. 
If $S \subseteq (0,1)$, then $I_S$ is the compact LOTS which
results by replacing each $x \in S$ by a pair of
neighboring points, $x^- < x^+$.
The \emph{double arrow space} is $I_{(0,1)}$.
\end{definition}

$I_S$ has no isolated points because
$0,1\notin S$.  The double arrow space is obtained by
splitting all points other than $0,1$.
$I_\emptyset = I$, and $I_{\QQQ \cap (0,1)}$ is homeomorphic
to the Cantor set.

\begin{lemma}
For each $S \subseteq (0,1)$, $I_S$ is a compact separable LOTS
with no isolated points.  $I_S$ is second countable iff $S$ is countable.
\end{lemma}

Now, let $Y = [0,1]$, let $S \subseteq (0,1)$, let $X = I_S$
and let $f : X \onto Y$ be the natural map.  Then
$f$ is $2$--tight by Lemma \ref{lemma-tight-LOTS}, but
$S = \{y \in Y : |f\iv\{y\}| \ge 2\}$ need not be countable, so
Theorems  \ref{thm-tight-metric},
\ref{thm-metric-embed}, and \ref{thm-metric-cover}
fail here when $S$ is uncountable (and hence $X$ is not metric).
However, one can apply these theorems in some generic extension,
to get a (perhaps strange) alternate proof that $f$ is $2$--tight.
Roughly, if $V[G]$ makes $S$ countable, then
$X,Y$ will both be compact metric in $V[G]$, so 
Theorem \ref{thm-tight-metric} implies that $f$ is 
$2$--tight in $V[G]$ (because $S$ is countable);
but then by absoluteness, $f$ is $2$--tight in $V$.
Absoluteness of tightness is discussed more precisely in Section \ref{sec-ab}.

The composition properties of tight maps are given by:

\begin{lemma}
\label{lemma-tight-comp}
Assume that $X,Y,Z$ are compact, $m,n$ are finite,
$f : X \onto Y$, and $g : Y \onto Z$.
Then:
\begin{itemizz}
\item[1.] If $g \circ f$ is $n$--tight then $g$ is $n$--tight.
\item[2.] If $f$ and $g$ are tight, then $g\circ f$ is tight.
\item[3.] If $f$ is $m +1$--tight and $g$ is $n +1$--tight,
then $g\circ f$ is $mn +1$--tight.
\end{itemizz}
\end{lemma}
\begin{proof}
(1) is trivial, and (2) is a special case of (3).

For (3), assume that $f$ is $m +1$--tight,  $g$ is $n +1$--tight,
and $g\circ f$ is not $mn +1$--tight; we shall derive a contradiction.
Fix disjoint closed $P_0,P_1, \ldots, P_{mn} \subseteq X$ with
$g(f(P_0)) \cap g(f(P_1)) \cap \cdots \cap g(f(P_{mn}))$ not scattered.
Shrinking $X,Y,Z$, and the $P_i$, we may assume WLOG that 
$X = P_0\cup P_1\cup \cdots \cup P_{mn}$ and that
$g(f(P_i)) = Z$ for each $i$.
For each $s \subseteq \{0,1,\ldots, mn\}$, let
$Q_s = \bigcap_{i\in s} f(P_i)$.
Shrinking the $P_i$, we may assume WLOG that each $Q_s \subseteq Y$
is either empty or not scattered; to see this, for a fixed $s$:
If $Q_s$ is scattered, then so is $g(Q_s)$; if $R$ is a perfect
subset of $Z \backslash g(Q_s)$, then we may replace $Z$ by $R$
and each $P_i$ by $P_i \cap f\iv(g\iv(R))$.

Now, using compactness of $ P_0\times P_1\times \cdots \times P_{mn}$,
as in the proof of Lemma \ref{lemma-tight-equiv}, fix $a_i \in P_i$
for $i \le mn$ such that
$g(f(\overline{U_0})) \cap\cdots\cap g(f(\overline{U_{mn}}))$ is not scattered
whenever each $U_i$ is a neighborhood of $a_i$.
Then at least one of the following two cases holds:

Case I.  Some $n+1$ of the  $f(a_0), \ldots, f(a_{mn})$ are different.
WLOG, these are $f(a_0), f(a_1), \ldots, f(a_{n})$.
Choose the $U_i$ so that the
$f(\overline{U_0}), f(\overline{U_1}), \ldots, f(\overline{U_{n}})$
are all disjoint.  But then 
$g(f(\overline{U_0})) \cap\cdots\cap g(f(\overline{U_{n}})) \supseteq
g(f(\overline{U_0})) \cap\cdots\cap g(f(\overline{U_{mn}}))$ is not scattered,
contradicting the $n +1$--tightness of $g$.

Case II.  Some $m+1$ of the  $f(a_0), \ldots, f(a_{mn})$ are the same.
WLOG, $f(a_0) = f(a_1) = \cdots = f(a_{m})$.
Let $s = \{0,1,\ldots, m\}$.
Then $Q_s \ne \emptyset$, so $Q_s = \bigcap_{i\le m} f(P_i)$ is not
scattered, contradicting the $m +1$--tightness of $f$.
\end{proof}

The ``$mn +1$'' in (3) cannot be reduced; for example, let
$Y = Z \times n$ and $X = Y \times m$, with $f,g$
the natural projection maps.

There is a similar result, with a similar proof, involving products:

\begin{lemma}
\label{lemma-tight-prod}
Assume that for $i = 0,1$:
$X_i, Y_i$ are compact,
$f_i:X_i \to Y_i$ is $(m_i + 1)$--tight,
$m_i \le n_i < \omega$, and
$|f_i\iv\{y\}| \le n_i$ for \emph{all} $y \in Y_i$.
Then $f_0 \times f_1 : X_0 \times X_1 \to Y_0 \times Y_1$ 
is $(\max(m_0 n_1 , m_1 n_0) + 1)$-- tight.
\end{lemma}
\begin{proof}
Let $L = \max(m_0 n_1 , m_1 n_0)$, and let $f = f_0 \times f_1$.
In view of Lemma \ref{lemma-tight-equiv}, it is sufficient
to fix any $L+1$ distinct points $a^0, a^1, \ldots, a^L \in X_0 \times X_1$
with all $f(a^\alpha) = b \in Y_0\times Y_1$,
and show that one can find neighborhoods $U^\alpha$ of $a^\alpha$ for 
$\alpha = 0, 1, \ldots , L$ such that $\bigcap_\alpha f(\overline{U^\alpha})$
is scattered.

Let $b = (b_0,b_1)$ and $a^\alpha = (a^\alpha_0, a^\alpha_1)$.

Note that although the $a^\alpha$ are all distinct points,
the $a^\alpha_0$ need not be all different and
the $a^\alpha_1$ need not be all different.
However,
$| \{a^\alpha_0 : 0 \le \alpha \le L \} | \ge m_0 + 1$:
If not, then using $f(a^\alpha) = b$ and
$|f_1\iv\{b_1\}| \le n_1$, we would have $L + 1 \le m_0 n_1$, a contradiction.
Likewise,  $| \{a^\alpha_1 : 0 \le \alpha \le L \} | \ge m_1 + 1$.

Now, using Lemma \ref{lemma-tight-equiv} and the fact that
each $f_i:X_i \to Y_i$ is $(m_i + 1)$--tight,
choose neighborhoods $U^\alpha_i$ of $a^\alpha_i$ such that 
$\bigcap_\alpha f(\overline{U^\alpha_i})$ is scattered for $i = 0,1$.
The $U^\alpha_i$ can depend just on the value of $a^\alpha_i$
(that is $a^\alpha_i =  a^\beta_i \to U^\alpha_i =  U^\beta_i$).
Finally, let $U^\alpha = U^\alpha_0 \times U^\alpha_1$
\end{proof}

The bound on the $|f_i\iv\{y\}|$ cannot be removed here.
For example, for each cardinal $\kappa$, one can find compact perfect
LOTSes $X_0,X_1,Y_0,Y_1$ with order-preserving $f_i : X_i \onto Y_i$
such that all point inverses have size at least $\kappa$.  Then
the $f_i$ are tight by Lemma \ref{lemma-tight-LOTS},
but $f_0 \times f_1$ is not $\kappa$--tight.

A variant of the product of maps is much simpler to analyze:

\begin{lemma}
\label{lemma-tight-prod-alt}
Assume that $\ell\in\omega$  and 
$f_i:X \to Y_i$ is $\kappa$--tight for each $i < \ell$,
where $X$ and the $Y_i$ are compact.
Then the map $x \mapsto (f_0(x), \ldots , f_{\ell - 1}(x) )$
from $X$ to $\prod_{i < \ell} Y_i$ is also $\kappa$--tight.
\end{lemma}

We now consider the opposite of tight maps:

\begin{definition}
\label{def-nowhere}
If $X,Y$ are compact and $f : X \to Y$, then
$f$ is \emph{nowhere tight} iff $f(X)$ is not scattered
and there is no closed $P \subseteq X$ such that
$f \res P$ is tight and $f(P)$ is not scattered.
\end{definition}

Note also that if $X,Y$ are metric compacta with $f : X \onto Y$
and $Y$ not scattered, then there is a Cantor set $P \subseteq X$
such that $f \res P$ is 1-1, so

\begin{lemma}
\label{lemma-metric-somewhere}
If $X,Y$ are compact and $f : X \to Y$ is nowhere tight,
then $X$ is not second countable. 
\end{lemma}

A further limitation on nowhere tight maps:

\begin{lemma}
\label{lemma-nwt}
If $f: X \to Y$ is nowhere tight, then $f$ is not weakly $\cccc$--tight.
\end{lemma}
\begin{proof}
We shall get a non-scattered $Q\subseteq Y$ and disjoint non-scattered sets
$P^k \subseteq X$ for  $k \in 2^\omega$ so that each $f(P^k) = Q$.
We shall build the $P^k$ and $Q$ by a tree argument.
Each $P^k$ will be non-scattered because it will be formed using a Cantor
tree of closed sets, so we shall actually get a doubly indexed
family.  So, we build
$Q_s \subseteq Y$ for $s \in 2^{< \omega}$ and 
$P_s^t \subseteq X$ for $s,t \in 2^{< \omega}$ with $\lh(s) = \lh(t)$
satisfying:
\begin{itemizz}
\item[1.] $P_s^t$ is closed, $f(P_s^t) =Q_s$, and $Q_s$ is not scattered.
\item[2.] The sets $Q_{s\cat0} , Q_{s\cat1}$ are disjoint subsets of $Q_s$.
\item[3.] The sets
$ P_{s\cat 0}^{t\cat 0} \, , \, P_{s\cat 0}^{t\cat 1} \, , \,
P_{s\cat 1}^{t\cat 0} \, , \, P_{s\cat 1}^{t\cat 1} $  are disjoint
subsets of  $P_s^t$.
\end{itemizz}
We construct these inductively.
$P_\one^\one$ and $Q_\one$ exist
(where $\one$ is the empty sequence) because $f(X)$ is not scattered.
Now, say we have $Q_s$ and $P_s^t$ for all $s,t$ with $\lh(s) = \lh(t) = n$.
Fix $s$.

First, get disjoint closed non-scattered
$\widetilde Q_{s\cat 0},\widetilde Q_{s\cat 1} \subseteq Q_s$, and let
$\widetilde P_{s\cat \mu}^{t} =
P_s^t \cap f\iv(\widetilde Q_{s \cat \mu})$ for
each $t$ of length $n$ and each $\mu = 0,1$.
Then, use ``nowhere tight'' $2^n$ times to get 
$ Q_{s \cat \mu} \subseteq \widetilde Q_{s \cat \mu}$ and
$P_{s\cat \mu}^{t\cat \nu} \subseteq \widetilde P_{s\cat \mu}^{t}$
for each $\mu,\nu = 0,1$ and each $t$ of length $n$ so that
each $f( P_{s\cat \mu}^{t\cat \nu} ) =  Q_{s \cat \mu}$
and each $Q_{s \cat \mu}$ is non-scattered.

For $h,k \in 2^\omega$, define
$Q_h = \bigcap_{n\in\omega} Q_{h\res n}$ and
$P_{h}^{k} = \bigcap_{n\in\omega} P_{h\res n}^{ k\res n}$,
let $Q = \bigcup\{Q_h : h \in 2^\omega\}$, and
let $P_h = \bigcup\{P_h^k  : k \in 2^\omega\}$ 
and $P^k = \bigcup\{P_h^k  : h \in 2^\omega\}$.
Then $f(P_h) = Q_h$ and $f(P^k) = Q$, and the
$\varphi$ of Definition \ref{def-weak-tight} sends
$P^k$ to $k \in 2^\omega$, with $\dom(\varphi) = \bigcup_k P_k$.
\end{proof}

\begin{corollary}
\label{cor-cantor}
If $X,Y$ are compact, $f : X \onto Y$,
$w(X) < \cccc$, $Y$ is metric and not scattered,
and $f$ is weakly $\cccc$--tight, then $X$ has a Cantor subset.
\end{corollary}
\begin{proof}
Since $f$ is not nowhere tight, we may assume, shrinking $X$ and $Y$,
that $f$ is tight.  Let $\kappa = w(X)$, and let $\BB$ be a base for $X$ with
$|\BB| = \kappa$.  Whenever $B_0,B_1 \in \BB$
with $\overline{B_0} \cap \overline{B_1} =\emptyset$,
let $S(B_0,B_1) = f(\overline{B_0} ) \cap f(\overline{B_1})$. 
Each $S(B_0,B_1)$ is scattered, and hence countable, 
so at most $\kappa$ points of $Y$ are in some $S(B_0,B_1)$,
so there is a $K \subseteq Y$ homeomorphic
to the Cantor set with $K$ is disjoint from all $S(B_0,B_1)$.
$|f\iv\{y\}| = 1$ for all $y \in K$, so 
$f\iv(K)$ is homeomorphic to $K$.
\end{proof}

Note that we have not yet given any examples of nowhere tight maps.
The argument of Corollary \ref{cor-cantor} shows that one
class of examples is given by:

\begin{example}
If $X,Y$ are compact, $f : X \onto Y$,
$w(X) < \cccc$, $Y$ is metric and not scattered,
and $X$ has no Cantor subset, then $f$ is nowhere tight.
\end{example}

Of course, under $\CH$, this class of examples is empty.
More generally, the class is empty 
under $\MA$ (or just the assumption that $\RRR$ is not the
union of $< \cccc$ meager sets), since then 
every non-scattered compactum of weight less than $\cccc$ contains
a Cantor subset (see \cite{HK3}).
However, by Dow and Fremlin \cite{DF},
it is consistent to have a non-scattered compactum $X$ of weight
$\aleph_1 < \cccc$ with no convergent $\omega$--sequences, and hence with no
Cantor subsets;  in the ground model, CH holds, and 
$X$ is any compact F-space (so $w(X)$ can be $\aleph_1$);
then, the extension adds any number of random reals.

A class of ZFC examples of nowhere tight maps with $w(X) = \cccc$ is
given by:

\begin{example}
If $X,Y$ are compact, $f : X \onto Y$,
$X$ is a compact F--space and
$Y$ is metric and not scattered,
then $f$ is nowhere tight.
\end{example}
\begin{proof}
Here, it is sufficient to prove that $f$ is not tight, since
any  $f\res P : P \onto f(P)$ will have the same properties.
Also, shrinking $Y$, we may assume that $Y$ has no isolated points.

First, choose a perfect $Q \subseteq Y$ which is nowhere dense in $Y$.
Then, choose a discrete set
$D = \{d_n : n \in \omega\} \subseteq Y \backslash Q$ with
$\overline D = D \cup Q$ and
each $f\iv\{d_n\}$ not a singleton.
Then, choose $x_n, z_n \in f\iv\{d_n\} $ with $x_n \ne z_n$.
Now, since $X$ is an F--space,
$\cl\{x_n : n \in \omega\}$ and
$\cl\{z_n : n \in \omega\}$ are two disjoint copies of $\beta\NNN$
in $X$ which map onto $\overline D$.
\end{proof}

\section{Dissipated Spaces}
\label{sec-dis}

Only a scattered compactum $X$ has the property that \emph{all} maps
from $X$ are tight: If $X$ is not scattered, then $X$ maps onto 
$[0,1]^2$; if we follow that map by the usual projection onto $[0,1]$,
we get a map from $X$ onto $[0,1]$ which is not
even weakly $\cccc$--tight.

The \emph{dissipated} compacta have the property that
\emph{unboundedly} many maps onto metric compacta are tight:

\begin{definition}
\label{def-fine}
Assume that $X,Y,Z$ are compact, $f : X \to Y$, and
$g : X \to Z$.  Then $f \le g$, or
$f$ is \emph{finer than} $g$, iff there is a
$\Gamma \in C(f(X),g(X))$ such that $g = \Gamma \circ f$.
\end{definition}

\begin{lemma}
\label{lemma-fine-equiv}
Assume that $X,Y,Z$ are compact, $f : X \to Y$, and
$g : X \to Z$.  Then $f \le g$ iff
$\; \forall x_1, x_2 \in X\, [ f(x_1) = f(x_2) \to g(x_1) = g(x_2)]$.
\end{lemma}
\begin{proof}
For $\leftarrow$,
let $\Gamma = \{( f(x),g(x)) : x \in X\} \subseteq f(X) \times g(X)$.
\end{proof}

\begin{definition}
\label{def-dis}
$X$ is \emph{$\kappa$--dissipated} iff $X$ is compact and
whenever $g : X \to Z$, with $Z$ metric,
there is a finer $\kappa$--tight $f : X \to Y$ for some metric $Y$.
$X$ is \emph{dissipated} iff $X$ is $2$--dissipated.
$X$ is \emph{weakly $\cccc$--dissipated} iff $X$ is compact and
whenever $g : X \to Z$, with $Z$ metric,
there is a finer weakly $\cccc$--tight $f : X \to Y$ for some metric $Y$.
\end{definition}

So, the $1$--dissipated compacta are the scattered compacta.
Metric compacta are trivially dissipated because we can take
$Y = X$, with $f$ the identity map.
Besides the spaces from \cite{Fed1, HK2, HK3}, an
easy example of a dissipated space is given by:

\begin{lemma}
\label{lemma-lots-dis}  
If $X$ is a compact LOTS, then $X$ is dissipated
\end{lemma}
\begin{proof}
Fix $g,Z$ as in Definition \ref{def-dis}.
On $X$, use $[x_1,x_2]$ for the closed interval
$[\min(x_1,x_2), \max(x_1,x_2)]$, and define $x_1 \sim x_2$
iff $g$ is constant on $[x_1,x_2]$.
Then $\sim$ is a closed equivalence relation, so define
$Y = X/\mosim$ with $f: X \onto Y$ the natural projection.
Then $Y$ is a LOTS and $f$ is order-preserving,
so $f$ is tight by Lemma \ref{lemma-tight-LOTS},
and $f \le g$ by Lemma \ref{lemma-fine-equiv}.
To see that $Y$ is metrizable, fix a metric on $Z$, and then, on $Y$,
define $d(f(x_1), f(x_2)) = \diam(g([x_1,x_2]))$.
\end{proof}

By Corollary \ref{cor-cantor}, if $w(X) < \cccc$ and $X$
is $\cccc$--dissipated and not scattered, then $X$ has a Cantor subset,
while the double arrow space is an example of an $X$ with
$w(X) = \cccc$ which is $2$--dissipated and has no Cantor subset.

Note that just having \emph{one} tight map $g$ from $X$ onto some
metric compactum $Z$ is not sufficient to prove that 
$X$ is dissipated, since the tightness of $g$ says nothing at all about
the complexity of a particular $g\iv\{z\}$.  Trivial counter-examples
are obtained with $|Z| = 1$ and $g$ a constant map.
However, if all $g\iv\{z\}$ are scattered, then just one tight $g$ is
enough:

\begin{lemma}
\label{lemma-one-enough}
Suppose that $g : X \to Z$ is $\kappa$--tight and all $g\iv\{z\}$
are scattered.  Fix $f : X \to Y$ with $f \le g$.
Then $f$ is $\kappa$--tight.  In particular, if $Z$ is also metric,
then $X$ is $\kappa$--dissipated.
\end{lemma}
\begin{proof}
Fix $\Gamma \in C(f(X),g(X))$ such that $g = \Gamma \circ f$.
Suppose that $\PP$ were a loose family for $f$ of size $\kappa$;
then we have $Q \subseteq f(X)$ with $Q = f(P)$ for all $P \in \PP$,
and $Q$ is not scattered.
But $\Gamma(Q)$ is scattered, since $g$ is $\kappa$--tight
and $g(P) = \Gamma(f(P)) = \Gamma(Q)$ for all $P \in \PP$.
It follows that we can fix $z\in Z$ with 
$Q \cap \Gamma\iv\{z\}$ not scattered.
But then $f(g\iv\{z\}) = \Gamma\iv\{z\}$ is not scattered,
which is impossible, since $g\iv\{z\}$ is scattered.
\end{proof}

We next consider the degree of dissipation of products:

\begin{lemma}
\label{lemma-prod-bad}
Let $X = A \times B$, where $A,B$ are compact,
$B$ is not scattered, and 
assume that for each $\varphi \in C(A, [0,1]^\omega)$
there is a $z \in [0,1]^\omega$
with $|\varphi\iv\{z\}| \ge \kappa$.  Then $X$ is not $\kappa$--dissipated.
If for each  $\varphi \in C(A, [0,1]^\omega)$ there is a $z$ such that
$\varphi\iv\{z\}$ is not scattered, then 
$X$ is not weakly $\cccc$--dissipated.
\end{lemma}
\begin{proof}
Since $B$ is not scattered, fix $h : B \onto [0,1]$, 
and define $g : X \onto [0,1]$ by $g(a,b) = h(b)$.
Now, fix any $f : X \to Y$  with $f$ finer than $g$ and
$Y$ metric.  We shall show that $f$ is not $\kappa$--tight.

Define $\widehat f : A \to C(B, Y)$ by
$(\widehat f(a))(b) = f(a,b)$.  Since the range of $\widehat f$ is compact
and hence embeddable in the Hilbert cube, we can fix 
$\zeta \in C(B,Y)$ such that $E := \{a : \widehat f(a) = \zeta\}$ has
size at least $\kappa$.
Let $Q = \zeta(B)$; $|Q| = \cccc$ by $f \le g$, so $Q$ is not scattered.
For $a \in E$, let $P_a = \{a\} \times B$.  Then $\{P_a : a \in E\}$
is a loose family of size at least $\kappa$.

The second assertion is proved similarly.
\end{proof}

Note that $A$ might be scattered; for example, $A$ could be the
ordinal $\kappa + 1$ (if $\kappa$ is uncountable and regular)
or the one point compactification of a discrete space of size
$\kappa$ (if $\kappa$ is uncountable).
$B$ may be second countable; for example $B$ can be the Cantor set.

A class of spaces $A$ to which Lemma \ref{lemma-prod-bad} applies
is produced by:

\begin{lemma}
\label{lemma-const}
Suppose that $f: \prod_{\alpha<\kappa} X_\alpha \to M$,
where $M$ is compact metric and, for each $\alpha$,
$X_\alpha$ is compact and not metrizable.  Then there are 
two-element sets $E_\alpha \subseteq X_\alpha$ for each $\alpha$ such that
$f$ is constant on $\prod_{\alpha<\kappa} E_\alpha $.
\end{lemma}
\begin{proof}
For $p \in \prod_{\alpha<\delta} X_\alpha $,
define $\widehat{f_p} : \prod_{\alpha\ge\delta} X_\alpha \to M$ by:
$\widehat{f_p}(q) = f(p \cat q)$.
Then inductively choose $E_\alpha$ so that
for all $\delta \le \kappa$, the functions $\widehat{f_p}$
are the same for all $p \in \prod_{\alpha<\delta} E_\alpha $.
Say $\delta < \kappa$ and we have chosen $E_\alpha$ for $\alpha < \delta$.
Let $g = \widehat{f_p}$ for
some (any) $p \in \prod_{\alpha<\delta} E_\alpha $, and define
$g^* \in C(X_\delta, C( \prod_{\alpha > \delta} X_\alpha, M)) $ by:
$(g^*(x))(q) = g(x \cat q)$.  Then $g^*$ maps $X_\delta$ into a metric
space of functions, so $\ran(g^*)$ is a compact metric space, so 
$g^*$ cannot be 1-1, so choose $E_\delta$ of size $2$
with $g^*$ constant on $E_\delta$.
\end{proof}

\begin{theorem}
\label{thm-bad-product}
Assume that each $X_k$ is compact:
\begin{itemizz}
\item[1.]
If $X_n$ is not scattered and $X_k$, for $k < n$, is not metrizable, then
$\prod_{k\le n}X_k$ is not $2^n$--dissipated.
\item[2.]
If each $X_k$ is not metrizable, then
$\prod_{k < \omega}X_k$ is not weakly $\cccc$--dissipated.
\end{itemizz}
\end{theorem}
\begin{proof}
For (1), apply Lemma \ref{lemma-prod-bad} with
$A = \prod_{k < n}X_k$  and $B = X_n$.
For (2), apply Lemma \ref{lemma-prod-bad} with
$A = \prod_{k < \omega}X_{2k}$  and
$B = \prod_{k < \omega}X_{2k + 1}$ .
\end{proof}

In (1), if all $X_k$ are scattered, then 
$\prod_{k\le n}X_k$ is scattered and hence dissipated.
As an example of (1) applied to LOTSes, if $S \subseteq (0,1)$
is uncountable, then
$(I_S)^2$ is not dissipated ($2$--dissipated),
$(I_S)^3$ is not $4$--dissipated, and
$(I_S)^4$ is not $8$--dissipated.
By Theorem \ref{thm-dis-prod}, these three spaces are, respectively,
$3$--dissipated, $5$--dissipated, and $9$--dissipated.
However, Lemma \ref{lemma-prod-bad}
shows that for any $\kappa$, we can
find a product of two LOTSes which is not $\kappa$--dissipated.

The following theorem will often suffice to compute the
degree of dissipation of a finite product of separable LOTSes:

\begin{theorem}
\label{thm-dis-prod}
Assume that $n$ is finite and $X_i$, for $i \le n$,
is a compact separable LOTS.
Then $\prod_{i\le n}X_i$ is $(2^n + 1)$--dissipated.
Furthermore, if all the $X_i$ are not scattered, and at most
one of the $X_i$ is second countable, then
$\prod_{i\le n}X_i$ is not $(2^n )$--dissipated.
\end{theorem}
\begin{proof}
Let $D_i \subseteq X_i$ be countable and dense.
Choose $f_i \in C(X_i, [0,1])$ such that $f_i$ is order-preserving
and is 1-1 on $D_i$ (such a function $f_i$ exists; see the proof
of Lemma 3.6 in \cite{HK1}).
Note that each $|f_i\iv\{y\}| \le 2$, and, by Lemma \ref{lemma-tight-LOTS},
each $f_i$ is $2$--tight.  Applying Lemma \ref{lemma-tight-prod}
and induction, $\prod_{i \le n}f_i$ is $(2^n + 1)$--tight.
Then $\prod_{i\le n}X_i$ is $(2^n + 1)$--dissipated by
Lemma \ref{lemma-one-enough}.

The ``furthermore'' is by Theorem \ref{thm-bad-product}.
\end{proof}

Next, we note that ``dissipated'' is a local property:

\begin{definition} 
Let $\KKKK$ be a class of compact spaces.
$\KKKK$ is \emph{closed-hereditary} iff every closed
subspace of a space in $\KKKK$ is also in $\KKKK$.
$\KKKK$ is \emph{local} iff $\KKKK$ is closed-hereditary \emph{and}
for every compact $X$:  if $X$ is covered by open sets whose
closures lie in $\KKKK$, then $X \in \KKKK$.
\end{definition} 

Classes of compacta which restrict cardinal functions
(first countable, second countable, countable tightness, etc.)
are clearly local, whereas the class of compacta which
are homeomorphic to a LOTS is closed-hereditary, but not local.
To prove that ``dissipated'' is local, we use as a preliminary lemma:

\begin{lemma}
\label{lemma-local-prelim}
Let $X$ be an arbitrary compact space, with $K \subseteq U \subseteq X$,
such that $U$ is open, $K$ is closed,
and $\overline U$ is $\kappa$--dissipated.
Fix $g : \overline U \to Z$, with $Z$ compact metric.
Then there is an $f : X \to Y$, with $Y$ compact metric,
$f$ $\kappa$--tight, and $f \res K \le g \res K$.
\end{lemma}
\begin{proof}
Fix $\varphi : X\to [0,1]$ with $\varphi(K) = \{0\}$ and
$\varphi(\partial U) = \{1\}$.
First get $f_0 : \overline U \to Y_0$, with $Y_0$ compact metric,
$f_0$ $\kappa$--tight, $f_0 \le g$, and $f_0 \le \varphi \res \overline U$
(just let $f_0$ refine $x \mapsto (g(x), \varphi(x))$).
Then $f_0(K) \cap f_0(\partial U) = \emptyset$.
Let $Y = Y_0 / f_0(\partial U)$, obtained by collapsing
$f_0(\partial U)$ to a point, $p$.
Let $f_1 :  \overline U \to Y$ be the natural map,
and extend $f_1$ to $f : X \to Y$ by letting $f_1(X\backslash U) = \{p\}$.
\end{proof}

\begin{lemma}
\label{lemma-local}
For any $\kappa$, the class of $\kappa$--dissipated compacta
is a local class.
\end{lemma}
\begin{proof}
For closed-hereditary:  Assume that $X$ is $\kappa$--dissipated
and $K$ is closed in $X$.
Fix $g : K \to Z$, with $Z$ metric.
Then we may assume that $Z \subseteq I^\omega$, so that 
$g$ extends to some $\widetilde g : X \to I^\omega$.
Then there is a $\kappa$--tight
$\widetilde f : X \to Y$ for some metric $Y$,
with $\widetilde f \le \widetilde g$.  If $f = \widetilde f \res K$,
then $f$ is  $\kappa$--tight and $f \le g$.

For local: Assume that $X = \bigcup_{i < \ell} U_i$,
where each $U_i$ is open and $\overline{U_i}$ is $\kappa$--dissipated.
Fix $g : X \to Z$, with $Z$ metric.
Choose closed $K_i \subseteq U_i$ such that $X = \bigcup_{i < \ell} K_i$.
Then apply Lemma \ref{lemma-local-prelim} and choose
$f_i : X \to Y_i$, with $Y_i$ compact metric,
$f_i$ $\kappa$--tight, and $f_i \res K_i \le g \res K_i$.
Then the map  $x \mapsto (f_0(x), \ldots , f_{\ell - 1}(x) )$ refines $g$,
and is  $\kappa$--tight by Lemma \ref{lemma-tight-prod-alt}.
\end{proof}

Many classes of compacta are closed under continuous images,
but this is not true in general of the
class of $\kappa$--dissipated spaces:

\begin{example}
\label{ex-three}
There is a continuous image of a $3$--dissipated space
which is not $\cccc$--dissipated.
\end{example}
\begin{proof}
Let $T = (D(\cccc) \cup \{\infty\}) \times 2^\omega$,
where $D(\cccc) \cup \{\infty\}$ is the 1-point compactification
of the ordinal $\cccc$ with the discrete topology.
Then $T$ is not $\cccc$--dissipated by Lemma \ref{lemma-prod-bad}.
Let $F_\alpha$, for $\alpha < \cccc$, be disjoint Cantor subsets
of $2^\omega$ such that for some $g: 2^\omega \onto 2^\omega$,
each $g(F_\alpha) = 2^\omega$.
Let $X = \{\infty\} \times 2^\omega
\ \cup\  \bigcup_{\alpha < \cccc} (\{\alpha\} \times F_\alpha)
\ \subseteq\  T$.  Then $X$ is $3$--dissipated by Lemma
\ref{lemma-one-enough} because the natural projection onto
$2^\omega$ is $3$--tight and all point inverses are
scattered (of size $\le 2$).
But also, $T$ is a continuous image of $X$ via the map
$\one \times g$, $(u,z) \mapsto (u,g(z))$.
\end{proof}

Of course, the continuous image of a $1$--dissipated (= scattered)
compactum is  $1$--dissipated.  We do not know about the
dissipated (= $2$--dissipated) spaces; perhaps $2$ is a special case.

\section{LOTS Dimension}
\label{sec-LOTS}
We shall apply the results of Section \ref{sec-dis} to products of LOTSes.
Each $I^n$ has dimension $n$ under any standard
notion of topological dimension, so that $I^{n+1}$ is
not embeddable into $I^n$.
Now, say we wish to prove such a result replacing $I$ by
some totally disconnected LOTS $X$.  Then standard
dimension theory gives all $X^n$ dimension $0$.
Furthermore, the result is false;
for example, $X^{n+1} \cong X^n$ if $X$ is the Cantor set.
However, if $X$ is the double arrow space, then 
$X^{n+1}$ is not embeddable into $X^n$.  To study this
further, we introduce a notion of LOTS dimension:

\begin{definition}
\label{def-lots-dim}
If $X$ is any Tychonov space, then $\Ldim_0(X)$
is the least $\kappa$ such that $X$ is embeddable into a product of
the form $\prod_{\alpha < \kappa} L_\alpha$, where each $L_\alpha$
is a LOTS.  Then $\Ldim(X)$, the \emph{LOTS dimension} of $X$,
is the least $\kappa$ such that every point in $X$
has a neighborhood $U$ such that $\Ldim_0(\overline U) \le \kappa$.
\end{definition}

\begin{lemma}
The class of compacta $X$ such that $\Ldim(X) \le \kappa$ is a local class.
\end{lemma}

If $X$ is any compact $n$--manifold, then $\Ldim(X) = n < \Ldim_0(X)$.
We follow the usual convention that the empty product
$\prod_{\alpha < 0} L_\alpha$ is a singleton,
so that $\Ldim(X) = 0$ iff $X$ is finite,
although $\Ldim_0(X) = 1$ if $1 < |X| < \aleph_0$.

\begin{lemma}
If $X$ is compact, infinite, and totally disconnected,
then $Ldim(X) = Ldim_0(X)$.
\end{lemma}
\begin{proof}
Use the fact that a disjoint sum of LOTSes is a LOTS.
\end{proof}

By Tychonov, $\Ldim(X) \le w(X)$, taking each $L_\alpha = I$.
In this section, we focus mainly on spaces whose LOTS dimension
is finite, although this cardinal function
might be of interest for other spaces.  For example,
$\Ldim(\beta\NNN) = 2^{\aleph_0}$; this is easily proved using
the theorem of Posp\'\i\v sil that there are points in 
$\beta\NNN$ of character $2^{\aleph_0}$.
We shall show (Lemma \ref{lemma-dim-product})
that $\Ldim((I_S)^n) = n$ whenever $S$ is uncountable.
When $S$ is countable, this is false if $S$ is dense in $I$
(then $(I_S)^n \cong I_S$ is the Cantor set) and true
if $S$ is not dense in $I$ (by standard dimension theory; not by
the results of this paper). 
More generally, we shall prove:

\begin{theorem}
\label{thm-ldim}
Let $Z_j$, for $1 \le j \le s$, be a compact LOTS.
Assume that $s  = r + m$, where $r,m \ge 0$.
For $r+1 \le j \le s$, assume that $Z_j$ has
either has an increasing or decreasing  $\omega_1$--sequence.
For $1 \le j \le r$, assume that there is a countable
$D_j \subseteq Z_j$ such that $\overline{D_j}$ is not scattered,
and assume that at most one of the $\overline{D_j}$ is
second countable.
Then $\Ldim(\prod_{j=1}^s Z_j) = s$.
\end{theorem}

The following lemma handles the case $r = s, m = 0$ if we
replace each $Z_j$ by $L_j = \overline{D_j}$.

\begin{lemma}
\label{lemma-dim-product}
Assume that $n$ is finite and $L_j$, for $j < n$,
is a compact separable LOTS.
Also, assume that all the $L_j$ are not scattered, and that at most
one of the $L_j$ is second countable.  Then $\Ldim(\prod_{j < n} L_j) = n$.
\end{lemma}
\begin{proof}
This is trivial if $n\le 1$, so assume that $n \ge 2$.
Clearly, $\Ldim(\prod_{j < n} L_j)  \le \Ldim_0(\prod_{j < n} L_j)  \le n$.
Also, by Theorem \ref{thm-dis-prod}, 
$\prod_{j < n} L_j$ is not $2^{n-1}$--dissipated.

To see that $\Ldim_0(\prod_{j < n} L_j)  \ge n$, assume
that we could embed $\prod_{j < n} L_j$ into
$\prod_{i < (n-1)}X_i$, where each $X_i$ is a LOTS.
Since the continuous image of a compact separable space
is compact and separable, we may assume that each $X_i$
is compact and separable, so that
$\prod_{i < (n-1)}X_i$ and $\prod_{j < n} L_j$,
are $(2^{n-2} + 1)$--dissipated by Theorem \ref{thm-dis-prod},
a contradiction since $2^{n-2} + 1 \le 2^{n-1}$.

Now, assume that $\Ldim(\prod_{j < n} L_j) <  n$.
Then we could cover $\prod_{j < n} L_j$ by finitely many open boxes,
each of the form $\prod_{j < n} U_j$, with each $U_j$ an open
interval in $L_j$, such that each open box satisfies
$\Ldim_0(\prod_{j < n} \overline{U_j}) <  n$.
But for at least one of these open boxes, the $\overline{U_j}$
would satisfy all the same hypotheses satisfied by the $L_j$,
so that we would again have a contradiction.
\end{proof}

In particular, if $L$ is the double arrow space, then
$L^{n+1}$ is not embeddable into $L^n$.  Similar results were obtained
by Burke and Lutzer \cite{BL} and Burke and Moore \cite{BM}
for the Sorgenfrey line $J$, which may be viewed as
$\{z^+ : z \in (0,1)\} \subseteq L$.  We do not see how to derive our results
directly from \cite{BL,BM}, since a map $\varphi : L^{n+1} \to L^n$
need not preserve order, so it does not directly yield a map
from $J^{n+1}$ to $J^n$.

We now extend Lemma \ref{lemma-dim-product} to include
LOTSes which have an
increasing or decreasing $\omega_1$--sequence.  First some preliminaries:

\begin{definition}
\label{def-seq}
$\seq{A}{n} = \{(\alpha_1, \ldots, \alpha_n) \in A^n :
\alpha_1 < \cdots < \alpha_n \}$, where
$1 \le n < \omega$ and $A \subseteq \omega_1$.
We give $\seq{A}{n} $ the topology it inherits from $(\omega_1)^n$.
The \emph{club filter} $\FF_n$ on $\seq{\omega_1}{n}$ is generated
by all the $\seq C n$ such that $C$ is club in $\omega_1$.
$\II_n$ is the dual ideal to $\FF_n$.
\end{definition}

\begin{lemma}
\label{lemma-clubs-decide}
If $B \subseteq \seq{\omega_1}{n}$ is a Borel set, then
$B \in \FF_n$ or $B \in \II_n$.
\end{lemma}
\begin{proof}
Since the $\II_n$ and $\FF_n$ are countably complete,
it is sufficient to prove this for closed sets $K$.
The case $n = 1$ is obvious, so we proceed by induction.
We assume the lemma for $n$, fix a closed
$K \subseteq  \seq{\omega_1}{(n+1)}$, and show that
$K \in \FF_{n+1}$ or $K \in \II_{n+1}$.
Applying the lemma for $n$:  For each $\alpha_0 < \omega_1$, 
choose $\nu(\alpha_0 ) \in \{0,1\}$ and a club 
$C_{\alpha_0} \subseteq (\alpha_0, \omega_1)$ such that for all
$(\alpha_1, \ldots, \alpha_n) \in \seq{ C_{\alpha_0}}{n}$:
\[
\nu(\alpha_0 ) = 0 \; \to\; (\alpha_0, \alpha_1, \ldots, \alpha_n) \notin K
\ \ ;  \ \   
\nu(\alpha_0 ) = 1 \; \to\; (\alpha_0, \alpha_1, \ldots, \alpha_n) \in K
\tag{$\ast$}
\]
Let $C = \{\delta : \delta \in \bigcap\{C_{\alpha_0} : \alpha_0 < \delta\}\}$.
Then $C$ is club and $(\ast)$ holds for all
$(\alpha_0, \alpha_1, \ldots, \alpha_n) \in \seq{C}{(n+1)}$.
Also, $D := \{\alpha_0 \in C: \nu(\alpha_0) = 1\}$ is closed
because $K$ is closed. 
$\seq{D}{(n+1)} \subseteq K$, so if $D$ is club, then $K \in \FF_{n+1}$.
If $D$ is bounded, then $C \backslash D$ contains a club,
and then $K \in \II_{n+1}$.
\end{proof}

\begin{definition}
If $L$ is a LOTS,
$f \in C( \seq{\omega_1}{m} , L)$, and $\psi \in C( \seq{\omega_1}{n} , L)$,
then  $\psi$ is \emph{derived from} $f$ iff
$n \ge m$ and for some $i_1, \ldots i_m$\textup:
$1 \le i_1< \cdots <i_m \le n$ and
$\psi( \alpha_1, \ldots, \alpha_n) = f(\alpha_{i_1}, \ldots, \alpha_{i_m})$
for all $(\alpha_1, \ldots, \alpha_n) \in \seq{\omega_1}{n}$.
Then a set $E \subseteq \seq{\omega_1}{n}$ is \emph{derived from} $f$ iff
$E$ is of the form
$\{\vec\alpha : \psi_1(\vec{\alpha}) < \psi_2(\vec{\alpha})\}$ or
$\{\vec\alpha : \psi_1(\vec{\alpha}) \le \psi_2(\vec{\alpha})\}$ or
$\{\vec\alpha : \psi_1(\vec{\alpha}) = \psi_2(\vec{\alpha})\}$,
where $\psi_1,\psi_2$ are derived from $f$.
\end{definition}

\begin{lemma}
\label{lemma-normal-form}
Suppose that $f \in C( \seq{\omega_1}{m} , L)$, where $L$ is a compact LOTS. 
Then there is a club $C$, a continuous $g:  C \to L$, and
a $j \in \{1,2,\ldots m\}$, such that for all
$\vec \alpha = (\alpha_1, \ldots, \alpha_m) \in \seq{C}{m}$, we
have $f(\vec\alpha) = g(\alpha_j)$, and $g$ is either 
strictly increasing or strictly decreasing or constant.
\end{lemma}
\begin{proof}
Applying Lemma \ref{lemma-clubs-decide}, and then restricting everything
to a club, we may make the following \emph{homogeneity} assumption:
for all $n \ge m$ and
all $E \subseteq \seq{\omega_1}{n}$ which are derived from $f$,
either $E = \emptyset$ or $E = \seq{\omega_1}{n}$.
Then, our club $C$ will be all of $\omega_1$.
We first consider the special cases $m = 1$ and $m = 2$.

For $m = 1$, we have $f \in C(\omega_1, L)$.  Applying homogeneity to
the three derived sets
$\{(\alpha,\beta) \in \seq{\omega_1}{2}  : f(\alpha) \circledast f(\beta)\}$,
where $\circledast$ is one of $<$, $>$, and $=$,
we see that $f$ is either
strictly increasing or strictly decreasing or constant.

Likewise, for $m > 1$, if we succeed in getting $f(\vec\alpha) = g(\alpha_j)$,
then $g$ must be either strictly increasing or strictly decreasing or constant.

Next, fix $f \in C( \seq{\omega_1}{2} , L)$.
If $\alpha < \beta < \gamma \to f(\alpha,\beta) = f(\alpha,\gamma)$,
then $f(\alpha,\beta) = g(\alpha)$, and we are done, so WLOG,
assume  $\alpha < \beta < \gamma \to f(\alpha,\beta) < f(\alpha,\gamma)$.
Let $B_\alpha = \{f(\alpha,\beta) : \alpha < \beta < \omega_1\}$,
which is a subset of $L$ of order type $\omega_1$.
Let $h(\alpha) = \sup(B_\alpha)$.  Fix $\alpha < \alpha' < \omega_1$.
There are now three cases; Cases II and III will lead to contradictions:

Case I.  $h(\alpha) = h(\alpha')$:  By continuity of $f$, there is
a club $C \subseteq (\alpha', \omega_1)$ such that
$f(\alpha,\beta) = f(\alpha',\beta) $ for all $\beta \in C$.
Applying homogeneity,
we have $\alpha < \alpha' < \beta \to f(\alpha,\beta) = f(\alpha',\beta) $,
so $f(\alpha,\beta) = g(\beta)$.

Case II.  $h(\alpha) < h(\alpha')$:  Fix $\beta$ such that
$\alpha < \alpha' < \beta$ and $f(\alpha',\beta) > f(\alpha,\gamma)$
for all $\gamma$.  Then by homogeneity,
$\alpha < \alpha' < \beta < \gamma \to f(\alpha,\gamma) < f(\alpha',\beta) $
for \emph{all} $\alpha , \alpha' , \beta , \gamma $.
Let $\alpha'$ be a limit and consider $\alpha \nearrow \alpha'$:  we get,
by continuity,
$\alpha' < \beta < \gamma \to f(\alpha',\gamma)\le f(\alpha',\beta) $,
contradicting $\alpha < \beta < \gamma \to f(\alpha,\beta) < f(\alpha,\gamma)$.

Case III.  $h(\alpha) > h(\alpha')$:  Fix $\beta$ such that
$\alpha < \alpha' < \beta$ and $f(\alpha,\beta) > f(\alpha',\gamma)$
for all $\gamma$.  Then by homogeneity,
$\alpha < \alpha' < \beta < \gamma \to f(\alpha',\gamma) < f(\alpha,\beta) $
for \emph{all} $\alpha , \alpha' , \beta , \gamma $.
Letting $\alpha \nearrow \alpha'$, we get a contradiction
as in Case II.

Finally, fix $m \ge 2$ and assume that the lemma holds for $m$.
We shall prove it for $m+1$, so fix $f \in C( \seq{\omega_1}{(m+1)} , L)$.
Temporarily fix
$(\alpha_1, \ldots, \alpha_{m-1}) \in \seq{\omega_1}{(m-1)}$,
and let $\widetilde f (\alpha_m, \alpha_{m+1}) =
f(\alpha_1, \ldots, \alpha_{m-1}, \alpha_m, \alpha_{m+1})$;
so $\widetilde f \in C( \seq{\,(\alpha_{m-1}, \omega_1)\, }{2} , L)$.  
Applying the $m=2$ case, $\widetilde f$ is really just a function of
one of its arguments, so that $f$ just depends on an 
$m$--tuple (either
$(\alpha_1, \ldots, \alpha_{m-1},  \alpha_{m+1})$
or
$(\alpha_1, \ldots, \alpha_{m-1}, \alpha_m)$), so we may now 
apply the lemma for $m$.
\end{proof}

It is easy to see from this lemma that $\Ldim( (\omega_1 + 1)^m) = m$,
but we now want to consider products of $(\omega_1 + 1)^m$ with
separable LOTSes.

\begin{lemma}
\label{lemma-normal-form-prod}
Suppose that $f \in C( X \times \seq{\omega_1}{m} , L)$,
where $L$ is a compact LOTS and $X$ is compact, nonempty,
first countable, and separable.
Then there is a club $C \subseteq \omega_1$,
a nonempty open $U \subseteq X$,
a $g \in C(\overline U \times C, L)$,
and a $j \in \{1,2,\ldots m\}$ such that
$f(x, \vec\alpha) = g(x, \alpha_j)$ for all
$\vec \alpha = (\alpha_1, \ldots, \alpha_m) \in \seq{C}{m}$
and all $x \in \overline U$, and such that either
\begin{itemizz}
\item[1.] for all $x \in \overline U$,
the map $\vec\alpha \mapsto f(x, \vec\alpha)$ is
constant on $\seq{C}{m}$ or
\item[2.] for all $x \in \overline U$,
the map $\xi \mapsto g(x,\xi)$ is strictly increasing on C, or
\item[3.] for all $x \in \overline U$,
the map $\xi \mapsto g(x,\xi)$ is strictly decreasing on C.
\end{itemizz}
\end{lemma}
\begin{proof}
First, let $K$ be the set of all $x$ such that 
$\vec\alpha \mapsto f(x, \vec\alpha)$
is constant on some set in $\FF_m$.  Then $K$ is closed,
since $X$ is first countable, so, replacing $X$ by some $\overline U$,
we may assume that $K = X$ or
$K = \emptyset$.  If $K = X$, then intersecting the clubs
for $x$ in a countable dense set, we get one club $C$ such
that (1) holds.

Now, assume that $K = \emptyset$.
Applying Lemma \ref{lemma-normal-form}, for each $x \in X$ choose
a club $C_x$, a $g_x \in C(C_x, L)$, and
$j_x \in \{1,2,\ldots m\}$ and a $\mu_x \in \{-1,1\}$
such that for all
$\vec \alpha = (\alpha_1, \ldots, \alpha_m) \in \seq{C_x}{m}$, we
have $f(x, \vec\alpha) = g_x(\alpha_{j_x})$, and each $g_x$ is either 
strictly increasing (when $\mu_x = 1$) or strictly decreasing
(when $\mu_x = -1$).

For each $j,\mu$, let
$H_j^\mu = \{x : j_x = j \ \&\ \mu_x = \mu \}$.  Then the $H_j^\mu$
are disjoint, and
they are also closed (since $K = \emptyset$).
Since $\bigcup_{j,\mu} H_j^\mu = X$, $U$ can be any nonempty $H_j^\mu$.
\end{proof}

In situations (2) or (3), we shall apply:

\begin{lemma}
\label{lemma-finite-range}
Suppose that $g \in C( X \times (\omega_1 + 1), L)$,
where $L$ is a compact LOTS and $X$ is compact, and suppose that
$g(x,\xi) < g(x,\eta)$ for each
$x \in X$ and each $\xi < \eta < \omega_1$.
Let $h(x) = g(x,\omega_1)$.  Then $h(X)$ is finite.
\end{lemma}
\begin{proof}
Assume that $h(X)$ is infinite.  Then, choose $c_n \in X$ for $n \in \omega$
such that the sequence  $\langle  h(c_n) : n \in \omega \rangle$
is either increasing strictly or decreasing strictly.
Let $c \in X$ be any limit point of $\langle c_n : n \in \omega \rangle$,
and note that $h(c_n) \to h(c)$.
Also note that $h(x) = \sup\{g(x,\xi) : \xi < \omega_1\}$ for every $x$.
Consider the two cases:

Case I.  $\langle  h(c_n) : n \in \omega \rangle$ is increasing strictly.
Then we can fix a large enough countable $\gamma$ such that
$g(c_n, \omega_1) < g(c_{n+1}, \gamma)$ for all $n$.  Then
we have the $\omega$--sequence,
$g(c_0, \gamma) < g(c_0, \omega_1)<
g(c_1, \gamma) < g(c_1, \omega_1) <
g(c_2, \gamma) < g(c_2, \omega_1)  < \cdots$,
whose limit must be $g(c, \gamma) = g(c, \omega_1)$, 
contradicting $g(c, \gamma) < g(c, \omega_1)$, 

Case II.  $\langle  h(c_n) : n \in \omega \rangle$ is decreasing strictly.
Then we can fix a large enough countable $\gamma$ such that
$g(c_n, \gamma) > g(c_{n+1}, \omega_1)$ for all $n$.  Then
we have the $\omega$--sequence,
$g(c_0, \omega_1) > g(c_0, \gamma)>
g(c_1, \omega_1) > g(c_1, \gamma) >
g(c_2, \omega_1) > g(c_2, \gamma)  > \cdots$,
whose limit must be $g(c, \omega_1) = g(c, \gamma)$, 
contradicting $g(c, \omega_1) > g(c, \gamma)$, 
\end{proof}

Now if $h(X)$ is finite, we can always shrink $X$ to a $\overline U$
on which $h$ is constant.  Then note that if $h(b) = h(c)$ and
$\xi \mapsto g(x,\xi)$ is always an increasing function,
then there is a club on which $g(b,\xi) = g(c,\xi)$.  Putting
these last two lemmas together, we get:

\begin{lemma}
\label{lemma-club}
Suppose that $f \in C( X \times (\omega_1 + 1)^m , L)$,
where $L$ is a compact LOTS and $X$ is compact, nonempty,
first countable, and separable.
Then there is a club $C \subseteq \omega_1$ and
a nonempty open $U \subseteq X$ such that either:
\begin{itemizz}
\item[1.] For some $j \in \{1,2,\ldots m\}$  and
some  continuous $g : C \to L$: 
$f(x, \vec\alpha) = g(\alpha_j)$
for all $x \in \overline U$ and all
$\vec \alpha  \in \seq{C}{m}$ and
$g$ is either 
strictly increasing or strictly decreasing, or \samepage
\item[2.] For some  $h \in C(\overline U , L)$:
$f(x, \vec\alpha) = h(x)$
for all $x \in \overline U$ and all
$\vec \alpha  \in \seq{C}{m}$.
\end{itemizz}
\end{lemma}

\begin{lemma}
\label{lemma-ldim-reduce}
Assume that $X$ is  compact, perfect, first countable, and separable,
and $\Ldim(X \times (\omega_1 + 1)^m) \le n$.  Then $n > m$
and there is a nonempty open $U \subseteq X$ such that
$\Ldim_0(\overline U) \le n-m$.
\end{lemma}
\begin{proof}
First, restricting everything to the closure of an open box,
we may assume that $\Ldim_0(X \times (\omega_1 + 1)^m) \le n$.

Fix a continuous 1-1 $f : X \times (\omega_1 + 1)^m \to \prod_{r = 1}^n L_r$,
where each $L_r$ is a compact LOTS.
Applying Lemma \ref{lemma-club} to the projections,
$f_r : X \times (\omega_1 + 1)^m \to  L_r$,
and permuting the $L_r$, we obtain a club $C$ and a $\overline U$
such that on $\overline U \times \seq{C}{m}$:
\[
f(x, \vec \alpha) = ( g_1(\alpha_{j_1}) , \ldots, g_p(\alpha_{j_p}) ,
h_1(x), \ldots, h_q(x) ) \ \ ,
\]
where $p + q = n$.  Then $\{j_1, \ldots, j_p\} = \{1, \ldots, m\}$,
since $f$ is 1-1.
Thus, $p \ge m$, so $q \le n-m$, and for any fixed $\vec \alpha$,
the map $x \mapsto (h_1(x), \ldots, h_q(x) )$ embeds
$\overline U$ into $\prod_{i = 1}^q L_{p +  i}$.
\end{proof}

\begin{proofof}{Theorem \ref{thm-ldim}}
Let $n = \Ldim(\prod_{j=1}^s Z_j)$.  Clearly $n \le s$.
To prove that $n \ge s$, we may
replace each $Z_j$ by a closed subset and assume that
$Z_j = \omega_1 + 1$ when $r+1 \le j \le s$, while
$Z_j = \overline{D_j}$ when $1 \le j \le r$. 
We may also assume that whenever $Z_j = \overline{D_j}$ is not second countable,
no open interval in $Z_j$ is second countable
(since there is always a closed subspace with this property).
Let $X = \prod_{j=1}^r Z_j$, and apply Lemma \ref{lemma-ldim-reduce}
to obtain $U \subseteq X$ with $\Ldim(\overline U) \le n-m$.
Since $\Ldim(\overline U) =  r $
by Lemma \ref{lemma-dim-product},
we have $r \le n-m$, so $s = r+m \le n$.
\end{proofof}

Note that this theorem does not cover all possible products of LOTSes.
For example, one can show by a direct argument that
$\Ldim( (\omega + 1) \times I_S) = 2$ whenever $S$ is uncountable, although
$(\omega + 1) \times I_S$ is dissipated, so the methods used in the proof
of Theorem \ref{thm-ldim} do not apply.
Also, Theorem \ref{thm-ldim} says nothing about Aronszajn lines,
which have neither an increasing or decreasing  $\omega_1$--sequence, nor
a countable subset whose closure is not second countable.
In particular, it is not clear whether one can have a product of three
compact Aronszajn lines which is embeddable into a product of two LOTSes.

In some sense, this ``dimension theory'' for products of
totally disconnected LOTSes
is more restrictive, not less restrictive, than the classical 
dimension theory for $I^n$, since there is also a limitation
on dimension-raising maps.  For example,
Peano \cite{PE} shows how to map $I$ onto $I^2$, but his map
has many changes of direction, so it does not define a map
from $I_S$ onto $(I_S)^2$.  In fact, this is impossible:

\begin{proposition}
If $S$ is uncountable, then there is no compact LOTS $L$
such that $L$ maps continuously onto $(I_S)^2$. 
\end{proposition} 
\begin{proof} 
Say $f : L \onto (I_S)^2$.  Replacing $L$ by a closed subset,
we may assume that $f$ is irreducible.  Then, $L$ must
be separable, since $(I_S)^2$ is separable.
It follows (see Lutzer and Bennett \cite{LB})
that $L$ is hereditarily separable,
which implies (by continuity of $f$)
that $(I_S)^2$ is hereditarily separable, which is false.
\end{proof}

We do not know whether, for example, one can map 
$L^2$ onto  $(I_S)^3$.  Again, we may assume that $L$ is separable,
so that $L^2$ is $3$--dissipated, while $(I_S)^3$ is not
even $4$--dissipated.  However, as we know from Example
\ref{ex-three}, a continuous image of a $3$--dissipated space
need not be even $\cccc$--dissipated.

\section{Measures, L-spaces, and S-spaces}
\label{sec-meas}

As usual, if $X$ is compact, a \emph{Radon measure} on $X$
is a finite positive regular Borel measure on $X$,
and if $f: X \to Y$ and $\mu$ is a measure on $X$,
then $\mu f\iv$ denotes the induced measure $\nu$ on $Y$,
defined by $\nu(B) = \mu(f\iv(B))$.
We shall prove some results relating $\mu$ to $\nu$ in the case that $f$ is
tight, and use this to prove that
Radon measures on dissipated spaces are separable.
We shall also make some remarks on compact L-spaces and S-spaces
which are dissipated.

\begin{definition}
For any space $X$, $\ro(X)$ denotes the \emph{regular open algebra} of $X$.
If $\BB$ is any boolean algebra and $b \in \BB$ with $b \ne \zero$,
then $b\down$ denotes the algebra $\{x \in \BB : x \le b\}$;
so $\one_{b\down} = b$.
A \emph{Suslin algebra} is an atomless ccc complete boolean algebra
which is $(\omega,\omega)$--distributive
\end{definition}

So, there is a Suslin tree iff there is a Suslin algebra.
We shall prove:

\begin{theorem}
\label{thm-suslin}
If $X$ is compact, ccc, not separable, and $\aleph_0$--dissipated,
then in $\ro(X)$ there is a non-zero $b$ such that
$b\down$ is a Suslin algebra.
\end{theorem}

Of course, this is well-known in the case where $X$ is a LOTS,
and is part of the proof that a Suslin line yields a Suslin tree.
Since a Suslin line is a compact L-space and is 
$2$--dissipated (by Lemma \ref{lemma-lots-dis}), we have

\begin{corollary}
There is an $\aleph_0$--dissipated compact L-space iff there is a Suslin line.
\end{corollary}

As usual, the \emph{support} of a Radon measure $\mu$ is the smallest closed
$H \subseteq X$ such that $\mu(H) = \mu(X)$.
For this $H$, $\ro(H)$ cannot be a Suslin algebra, so

\begin{corollary}
If $X$ is $\aleph_0$--dissipated, then
the support of every Radon measure on $X$ is a separable
topological space.
\end{corollary}

In these two corollaries, the ``$\aleph_0$'' cannot be
replaced by ``$\aleph_1$'', since
the usual compact L-space construction shows the following
(see Section \ref{sec-inv} for a proof):

\begin{proposition}
\label{prop-L}
$\CH$ implies that there is a compact L-space $X$ which is both
$\cccc$--dissi\-pated and the support of a Radon measure $\mu$.
Furthermore, $\mu$ is atomless, and, in $X$, the ideals
of null subsets, meager subsets, and separable subsets
all coincide.
\end{proposition}

Turning to compact S-spaces, the usual CH construction \cite{JKR} yields
one which is scattered, and hence dissipated.  Less trivially,
the construction of Fedor\v cuk \cite{Fed1} shows, under $\diamondsuit$,
that there is a dissipated compact S-space with no isolated points 
and no non-trivial convergent $\omega$--sequences;
see Section \ref{sec-inv} for further remarks on this construction.

\begin{proofof}{Theorem \ref{thm-suslin}}
Since $X$ is ccc, we may replace $X$ by some regular closed set
and assume that $X$ is nowhere separable --- that is, the closure
of every countable subset is nowhere dense.
Assume that in $\ro(X)$ no  $b\down$ is Suslin, and we shall
derive a contradiction.

Since $X$ is ccc, the fact that no $b\down$ is Suslin implies
that there are open $F_\sigma$ sets $V_n^j$ for $n,j \in \omega$ such that
for each $n$, the $V_n^j$ for $j \in \omega$ are disjoint and
$\bigcup_j V_n^j$ is dense, and such that for each $\varphi \in \omega^\omega$,
$\bigcap_n V_n^{\varphi(n)}$ has empty interior.
There is then a compact metric $Y$ and an $f : X \onto Y$
such that $V_n^j  =  f\iv(f(V_n^j))$ for each $n,j$.
Note that this implies that each $f(V_n^j)$ is open,
since $f(V_n^j)  = Y \setminus f(X \backslash V_n^j)$.

Replacing $f$ by a finer map, we may also assume that $f$
is $\aleph_0$--tight.

Observe that $f\iv\{y\}$ is nowhere dense for each $y \in Y$,
since either $f\iv\{y\} \subseteq \bigcap_n V_n^{\varphi(n)}$
for some $\varphi \in \omega^\omega$, or
$f\iv\{y\} \subseteq X \setminus \bigcup_j V^j_n$ for some $n$.

Now, construct open $U_s \subseteq X$ and closed $K_s \subseteq X$
for $s \in 2^{<\omega}$ as follows:
$U_{(\,)} = X$, and each $\overline{U_{s\cat i}} \subseteq U_s \backslash K_s$,
with $f( \overline{U_{s\cat 0}}) \cap f( \overline{U_{s\cat 1}}) = \emptyset$.
Also, $K_s \subseteq \overline{U_s}$, with
$f(K_s) = f(  \overline{U_s}) $ and $f \res K_s \onto f(K_s)$ irreducible.
Note that $K_s$ is separable and $\overline{U_s}$ is nowhere separable,
so that the construction can continue.
More specifically, to choose $U_{s\cat 0}$ and $U_{s\cat 1}$:
First, find $p_0, p_1 \in U_s \backslash K_s$ such that $f(p_0) \ne f(p_1)$;
this is possible since otherwise we would have
$f(U_s \backslash K_s) \subseteq \{y\}$, contradicting the
fact that $f\iv\{y\}$ is nowhere dense.
Next, find open $W_i \subseteq Y$ with $f(p_i) \in W_i$ and
$\overline{W_0} \cap \overline{W_1} = \emptyset$.
Then, choose $U_{s\cat i}$ with
$U_{s\cat i}\subseteq \overline{U_{s\cat i}} \subseteq
(U_s\backslash K_s) \cap f\iv(W_i)$.

Let $Q_n = \bigcup\{f(K_s) : s \in 2^n\}$, and let
$Q = \bigcap_n Q_n$, which is a non-scattered subset of $Y$.
Let $P_n =  f\iv(Q) \cap \bigcup\{K_s : s \in 2^n\}$.
Then the $P_n$ are disjoint and each $f(P_n) = Q$,
contradicting the $\aleph_0$--tightness of $f$.
\end{proofof}

To study measures further, we use the following standard definitions:

\begin{definition}
If $\mu$ is any finite measure on $X$, then $\ma(\mu)$ denotes the
\emph{measure algebra} of $\mu$ --- that is, the algebra
of measurable sets modulo the null sets.
If $f: X \to Y$, 
$\mu$ is a finite measure on $X$, and $\nu = \mu f\iv$, then
$f^*: \ma(\nu) \to \ma(\mu)$ is defined by $f^*([A]) = [f\iv(A)]$. 
\end{definition}

$\ma(\mu)$  is a complete metric
space with metric $d([A], [B]) = \mu(A \Delta B)$,
where $[A],[B]$ denote the equivalence classes of the sets $A,B$.
Note that we do not require $f$ to be onto here, although
$Y \backslash f(X)$ is a $\nu$--null set.
$f^*$ is an isometric isomorphism
onto some complete subalgebra $ f^*(\ma(\nu)) \subseteq \ma(\mu)$.

As usual, a measure $\mu$ on $X$ is \emph{separable} iff
$L^p(\mu)$ is a separable metric space for some (equivalently, for all)
$p \in [1,\infty)$.  Also $\mu$ is separable iff 
$\ma(\mu)$ is a separable metric space iff
$\ma(\mu)$ is countably generated as a complete boolean algebra.
Separability of $\mu$ is not related in any simple way to 
the separability of any topology that $X$ may have.
Following \cite{DK2}:

\begin{definition}
\label{def-MS}
$\MS$ is the class of all compact spaces $X$ such that every
Radon measure on $X$ is separable.
\end{definition}

We shall prove:

\begin{theorem}
\label{thm-dis-ms}
If $X$ is a weakly $\cccc$--dissipated space then $X$ is in $\MS$.
\end{theorem}

In view of Lemma \ref{lemma-lots-dis}, Theorem \ref{thm-dis-ms} generalizes
the result from \cite{DK2} that every compact LOTS is in $\MS$.
Note that a space in $\MS$ need not be $\cccc$--dissipated.  For example,
$\MS$ is closed under countable products (see \cite{DK2}), but an
infinite product of non-metric compacta is never
weakly $\cccc$--dissipated (see Theorem \ref{thm-bad-product}).

Theorem \ref{thm-dis-ms} will be an easy corollary
of some general results about measures induced by weakly
$\cccc$--tight $f: X \onto Y$, where $X,Y$ are compact.
Say $\mu$ is a Radon measure on $X$, with $\nu = \mu f\iv$.
Even if $f$ is tight (i.e., $2$--tight), the separability of 
$\nu$ does not imply the separability of $\mu$;
for example, $\nu$ may be a point mass concentrating on $\{y\}$,
in which case $\mu$ can be any measure supported on $f\iv\{y\}$ with
$\mu(f\iv\{y\}) = \nu\{y\}$.
However, if $\nu$ is atomless, then
the form of $\nu$ will restrict the form of $\mu$.
There are really two kinds of ways that $\nu$ might determine $\mu$.
We shall denote the stronger way as ``$X$ is skinny'' and the
weaker way as ``$X$ is slim''.  We shall define
``skinny'' and ``slim'' also for arbitrary closed subsets of $X$:

\begin{definition}
\label{def-skinny -slim}
Suppose that $X,Y$ are compact,
$f: X \to Y$, $\mu$ is a  Radon measure on $X$, and $\nu = \mu f\iv$. Then:
\begin{itemizn}{43}
\item $X$ is \emph{skinny} with respect to $f,\mu$ iff 
for all closed $K \subseteq X$, $\mu(K) = \nu(f(K))$.
\item $X$ is \emph{slim} with respect to $f,\mu$ iff 
$f^*: \ma(\nu)  \to \ma(\mu)$
maps \emph{onto} $\ma(\mu) $.
\end{itemizn}
If $H$ is a closed subset of $X$, then we say that $H$ is
\emph{skinny} \textup(resp., \emph{slim}\textup) with respect to $f,\mu$ iff 
$H$ is skinny \textup(resp., slim\textup) with respect to
$f\res H,\mu\res H$.
\end{definition}

Note that the equation $\mu(K) = \nu(f(K))$ shows that if $X$
is skinny, then $\nu$ determines $\mu$; there is no Radon measure
$\mu' \ne \mu$ such that $\nu = \mu' f\iv$.

\begin{lemma}
\label{lemma-skinny-to-slim}
If $X$ is skinny with respect to  $f,\mu$, then $X$ is slim.
\end{lemma}
\begin{proof}
If $K  \subseteq X$ is closed, then 
$\mu(K) = \mu(f\iv(f(K)))$ implies
that $[K]  = [f\iv(f(K))] =f^*([f(K)])$ in $\ma(\mu)$. 
Thus, $[K] \in \ran( f^* )$ for all
closed $K  \subseteq X$, which implies that $f^* $ is onto.
\end{proof}

The converse is false.  For example, suppose that $H$ is a closed
subset of $X$ such that $\mu$ is supported on $H$ and $f \res H$
is 1-1.  Then $X$ is slim, since $\ma(\mu) \cong \ma(\mu \res H)$,
but $X$ need not be skinny, since there may well be closed $K$
disjoint from $H$ with $X = f\iv(f(K))$; then $\mu(K) = 0$
but $\nu(f(K)) = \mu(X)$.
In this example, $H$ is skinny with respect to $f, \mu$.  Some examples
of skinny sets on which the function $f$ is not 1-1 are given by:

\begin{lemma}
\label{lemma-tight-to-skinny}
Suppose that $X,Y$ are compact,
$f: X \to Y$ is tight,
$\mu$ is a Radon measure on $X$, and $\nu = \mu f\iv$ is atomless. 
Then $X$ is skinny with respect to $f,\mu$.
\end{lemma}
\begin{proof}
If $X$ is not skinny, fix a 
closed $K \subseteq X$ with $\mu(K) < \nu(f(K))$, 
so that $\mu( f\iv(f(K)) \setminus K) > 0$.
Then choose a closed
$L \subseteq  f\iv(f(K)) \setminus K$ with $\mu(L) > 0$.
Then $K,L$ are disjoint in $X$ and
$\nu( f(K) \cap f(L)) = \nu(f(L)) \ge \mu(L) > 0$,
so $ f(K) \cap f(L)$ cannot be scattered, since $\nu$ is atomless,
so $f$ is not tight.
\end{proof}

One cannot replace ``tight'' by ``$3$--tight'' here.
For example, say $X = Y \times \{0,1\}$, with $f$ the natural
projection, which is $2$--tight.  If $\nu$ is any Radon
measure on $Y$,
and on $X$ we let $\mu( E_0 \times \{0\} \;\cup\;  E_1 \times \{1\}) =
\frac12 (\nu(E_0) + \nu(E_1))  $, then $X$ is not skinny (or even slim).
Here, $X$ is the union of two skinny subsets,
and this situation generalized to:

\begin{lemma}
\label{lemma-tight-skinny}
Suppose that $X,Y$ are compact,
$f: X \onto Y$ is $\aleph_0$--tight 
and $\mu$ is a Radon measure on $X$ with $\mu f\iv$ atomless.
Then there is a countable family $\HH$ of disjoint skinny subsets of
$X$ such that $\mu(X) = \sum\{\mu(H) : H \in \HH\}$.
\end{lemma}
\begin{proof}
If this fails, then the usual exhaustion argument lets us shrink
$X$ and assume that $\mu(X) > 0$ and there are no closed skinny $H \subseteq X$
of positive measure.  We now build an infinite loose family as follows:

Construct a tree of closed $H_s \subseteq X$ for $s \in 2^{< \omega}$;
so $H_{s\cat 0}, H_{s\cat 1}$ will be disjoint closed subsets of $H_s$,
and also $f(H_{s\cat 0}) \cap f(H_{s\cat 1}) = \emptyset$.
Each $H_s$ will have positive measure.  $H_{(\,)}$ can be $X$.

Given $H_s$:  Since $H_s$ is not skinny, we can choose a closed
$K_s \subset H_s$ with $\mu(H_s \cap (f\iv(f(K_s)) \setminus K_s)) > 0$.
Then, since $\mu$ is regular and $\mu f\iv$ is atomless,
we can choose closed
$H_{s\cat 0}, H_{s\cat 1} \subseteq H_s \cap (f\iv(f(K_s)) \setminus K_s)$
with each $\mu(H_{s\cat i}) > 0$ and
$f(H_{s\cat 0}) \cap f(H_{s\cat 1}) = \emptyset$.

Now, let $Q_n = \bigcup\{f(H_s) : s \in 2^n\}$ and let
$Q = \bigcap_n Q_n$; so, $Q$ is non-scattered.  Let
$P_n = f\iv(Q) \cap \bigcup\{K_s : s \in 2^n\}$.
Then $\{P_n : n \in \omega\}$ is a loose family.
\end{proof}

It follows that the measure algebra of $\mu$ is a countable sum 
of measure algebras isomorphic to algebras derived from measures on $Y$.
Note that the $K_s$ in this proof may be null sets, so
one cannot split them also to obtain a loose family of size $\cccc$,
as we did in the proof of Lemma \ref{lemma-nwt}.
In fact, the L-space of Proposition \ref{prop-L} shows that one cannot
weaken ``$\aleph_0$--tight'' to ``$\aleph_1$--tight'' in this lemma.
To see this, note that $\mu$ is a separable measure on $X$
by Theorem \ref{thm-dis-ms}, so one can get an $f: X \onto Y$ such
that $Y$ is compact metric, $\nu = \mu f\iv$ atomless,
and $f^*(\ma(\nu)) = \ma(\mu)$.
Since  $X$ is $\aleph_1$--dissipated, one can refine $f$
and assume also that $f$ is $\aleph_1$--tight.
Now, if $H$ is skinny, let $K$ be a closed subset of $H$ such that
$f(K) = f(H)$ and $f\res K : K \onto f(H)$ is irreducible.
Then $K$ is separable and hence null (by the properties of $X$),
and $\mu(H) = \mu(K)$ (since $H$ is skinny), so $\mu(H) = 0$.
Thus, there cannot be a family $\HH$ as in Lemma \ref{lemma-tight-skinny}.

However, the analogous result with ``slim''
(Theorem \ref{thm-tight-to-simple}) just uses $\cccc$--tightness.

\begin{definition}
Suppose that $X,Y$ are compact,
$f: X \to Y$, and  $\mu$ is a Radon measure on $X$.
Then $X$ is \emph{simple} with respect to $f,\mu$
iff there is a countable disjoint family
$\HH$ of slim subsets of $X$ such that $\sum\{\mu(H) : H \in \HH\} = \mu(X)$.
\end{definition}

We shall prove:

\begin{theorem}
\label{thm-tight-to-simple}
Suppose that $X,Y$ are compact,
$f: X \to Y$, and $\mu$ is a Radon measure on $X$, with $\nu = \mu f\iv$,
and suppose that $X$ is not simple with respect to $f,\mu$.
Then there is a $\varphi : \dom(\varphi) \to 2^\omega$,
where $\dom(\varphi)$ is closed in $X$,
such that for some closed $Q \subseteq Y$,
$\nu(Q) > 0$ and $\varphi(f\iv\{y\}) = 2^\omega$ for all $y \in Q$.
In particular, if $\nu$ is atomless, then $f$ is not weakly $\cccc$--tight.
\end{theorem}

In proving this, the notion of  \emph{conditional expectation}
(see \cite{Halmos}, \S48)
will be useful in comparing the induced measure
$(\mu \res S) f\iv$ to $\nu$ for various $S \subseteq X$:

\begin{definition}  
Suppose that $f: X \to Y$, with $X,Y$ compact,
$\mu$ is a measure on $X$ and $\nu = \mu f\iv$.
If $S$ is a measurable subset of $X$, then 
the \emph{conditional expectation},
$\EEE(S | f) = \EEE_\mu(S | f)$, is the measurable
$\varphi: Y \to [0,1]$ defined so that
$\int_A \varphi(y) \, d\nu(y) = \mu( f\iv(A) \cap S)$ for all measurable
$A \subseteq Y$.
\end{definition}

Of course, $\varphi$ is only defined up to equivalence in $L^\infty(\nu)$.
Conditional expectations are usually defined for probability measures,
but they make sense in general for finite measures;
actually, $\EEE_\mu(S | f) = \EEE_{c\mu}(S | f)$ for any non-zero $c$.
Note that
$\int_A \varphi(y) \, d\nu(y) = \int_{f\iv(A)} \varphi(f(x)) \, d\mu(x)$.
We may also characterize $\varphi = \EEE_\mu(S | f)$ by the equation:
\[
\int_S g(f(x)) \, d\mu(x) = \int_X \varphi(f(x))\, g(f(x))  \, d\mu(x) =
\int_Y \varphi(y)\, g(y)  \, d\nu(y)  \ \ .
\]
for $g \in L^1(Y, \nu)$.  $\varphi$ is obtained either by the Radon-Nikodym
Theorem, or, equivalently,
by identifying $(L^1(Y, \nu))^* $ with $L^\infty(Y, \nu)$,
since $\Gamma(g) := \int_S g(f(x))\, dx$ defines $\Gamma \in (L^1(Y, \nu))^* $,
with $\|\Gamma \| \le 1$.

Now, given $\mu$ on $X$ and $f : X \to Y$, we shall consider
various closed subsets $H \subseteq X$ while studying the tightness properties
of $f$.  When $S \subseteq H \subseteq X$,
one must be careful to distinguish $\EEE_\mu(S|f)$
(computed using $\mu$ and $f : X \to Y$)
from $\EEE_{\mu \res H} (S \mid f\res H)$ 
(computed using $\mu\res H$ and $f\res H : H \to Y$).
These are related by:

\begin{lemma}
\label{lemma-realate-exp}
Suppose that $f: X \to Y$, with $X,Y$ compact,
$H$ is a closed subset of $X$, and
$\mu$ is a Radon measure on $X$.  Let $S$ be a measurable subset of $H$.
Then  $\EEE_\mu(S|f) = \EEE_\mu(H|f) \cdot \EEE_{\mu \res H} (S \mid f\res H)$.
\end{lemma}
\begin{proof}
Let $\varphi = \EEE_\mu(S|f)$, $\psi =  \EEE_\mu(H|f)$,
and $\gamma = \EEE_{\mu \res H} (S \mid f\res H)$.
We may take these to be bounded Borel-measurable functions from $Y$ to $\RRR$.
For any bounded Borel-measurable $g: Y \to \RRR$, we have
\begin{align*}
& \int_S g(f(x)) \, d\mu(x) = \int_X \varphi(f(x))\, g(f(x))  \, d\mu(x) \\
& \int_H g(f(x)) \, d\mu(x) = \int_X \psi(f(x))\, g(f(x))  \, d\mu(x) \\
& \int_S g(f(x)) \, d\mu(x) = \int_H \gamma(f(x))\, g(f(x))  \, d\mu(x)  =
  \int_X \psi(f(x))\, \gamma(f(x))\, g(f(x)) \, d\mu(x),
\end{align*}
which yields $\varphi = \psi \gamma$.
\end{proof}

We now relate conditional expectations to slimness:

\begin{lemma}
\label{lemma-ran-fstar}
Suppose that $X,Y$ are compact,
$f: X \to Y$, and $\mu$ is a measure on $X$, with $\nu = \mu f\iv$.
Let $S \subseteq X$ be measurable.  Then $[S] \in \ran(f^*)$ iff
$[\EEE(S|f)] = [\cchi_T]$ for some measurable $T \subseteq Y$,
in which case $[S] = f^*([T])$.
\end{lemma}
\begin{proof}
For $\rightarrow$:  If $[S] = f^*([T])$ then $\mu(S \Delta f\iv(T)) = 0$,
which implies $\EEE(S | f) = \EEE(f\iv(T) | f) = \cchi_T$.

For $\leftarrow$:  If $[\EEE(S|f)] = [\cchi_T]$ then
$ \mu( f\iv(A) \cap S) = \nu(A \cap T)$ for all measurable $A \subseteq Y$.
Setting $A = Y \backslash T$, we get $ \mu(S \backslash  f\iv(T)) = 0$,
so $[S] \le [f\iv(T)]$. Setting $A = T$, we get
$\mu(S \cap f\iv(T)) =  \nu(T) = \mu(f\iv (T))$, so 
$[S] \ge [f\iv(T)]$. 
\end{proof}

In particular, $X$ is slim with respect to $f,\mu$
iff every $\EEE(S | f)$ is the characteristic function of a set;
this remark will be useful when applied also to $\mu\res H$
for various $H \subseteq X$.  

\begin{lemma}
\label{lemma-non-slim-split}
Suppose that $X,Y$ are compact,
$f: X \to Y$, and $\mu$ is a measure on $X$, with $\nu = \mu f\iv$,
and suppose that $X$ is not slim with respect to $f,\mu$.
Then there are disjoint closed $H_0, H_1 \subseteq X$
with $f(H_0) = f(H_1) = K$, such that $\nu(K) > 0$ and, for $i = 0,1$,
$\; 0 < \EEE(H_i | f)(y) < 1$ for a.e.\@ $y \in K$.
\end{lemma}
\begin{proof}
First,
let $\widetilde H_0 \subseteq X$ be closed with $[H_0] \notin \ran(f^*)$.
We can then, by Lemma \ref{lemma-ran-fstar}, get a closed
$\widetilde K \subseteq f( \widetilde H_0)$  with $\nu(\widetilde K) > 0$ and
$\EEE(\widetilde H_0 | f)(y) \in (0,1)$ for a.e.\@ $y \in \widetilde K$. 
Then, choose a closed
$\widetilde H_1 \subseteq f\iv(\widetilde K) \backslash \widetilde H_0$
with $\mu( \widetilde H_1 ) > 0$.
Then, choose a closed $K \subseteq \widetilde f( \widetilde H_1)$
with $\nu(K) > 0$ and $\EEE(\widetilde H_1 | f)(y) > 0$ for a.e.\@ $y \in K$,
and let $H_i = \widetilde H_i \cap f\iv(K)$.
\end{proof}

We now consider the opposite of slim:

\begin{definition}
$X$ is \emph{nowhere slim} with respect to $f, \mu$ iff there is no
closed $H \subseteq X$ with $\mu(H) > 0$
such that $H$ is slim with respect to $f, \mu$.
\end{definition}

\begin{lemma}
\label{lemma-nowhere-slim-split}
Suppose that $X,Y$ are compact,
$f: X \to Y$, and $\mu$ is a measure on $X$, with $\nu = \mu f\iv$,
and suppose that $X$ is nowhere slim with respect to $f,\mu$.
Fix $\varepsilon > 0$.
Then there are disjoint closed $H_0, H_1 \subseteq X$
with $f(H_0) = f(H_1) = K$, such that $\nu(Y \backslash K) <  \varepsilon$
and, for $i = 0,1$,
$\; 0 < \EEE(H_i | f)(y) < 1$ for a.e.\@ $y \in K$.
\end{lemma}
\begin{proof}
Fix $\KK$ such that
\begin{itemizz}
\item[1.] $\KK$ is a disjoint family of non-null closed subsets of $Y$.
\item[2.] For $K \in \KK$, there are disjoint closed $H^K_0, H^K_1 \subseteq X$
with $f(H^K_0) = f(H^K_1) = K$,
and, for $i = 0,1$,
$\; 0 < \EEE(H^K_i | f)(y) < 1$ for a.e.\@ $y \in K$.
\item[3.] $\KK$ is maximal with respect to (1)(2).
\end{itemizz}
Then $\KK$ is countable.  If $\nu(Y \backslash \bigcup \KK) = 0$,
choose a finite $\KK' \subseteq \KK$ such that
$\nu(Y \backslash \bigcup \KK') < \varepsilon$, set
$K = \bigcup \KK'$, and set $H_i = \bigcup \{H_i^K : K \in \KK'\}$.
If $\nu(Y \backslash \bigcup \KK) \ne 0$, choose a closed 
$E \subseteq Y \backslash \bigcup \KK $ with $\nu(E) > 0$,
and use Lemma \ref{lemma-non-slim-split} to 
derive a contradiction from maximality of $\KK$ and the
fact that $f\iv(E)$ is not slim.
\end{proof}

We can now use a tree argument to prove Theorem \ref{thm-tight-to-simple}:

\begin{proofof}{Theorem \ref{thm-tight-to-simple}}
Since $f$ is not simple, there must be a closed $H \subseteq X$
such that $H$ is nowhere slim with respect to $\mu \res H, f \res H$.
Restricting everything to $H$, we may assume that $X$ itself is
nowhere slim.  Also, WLOG $\mu(X) = \nu(Y) = 1$ and $f(X) = Y$.
Now, get $P_s \subseteq X$ for $s \in 2^{<\omega}$ and
$Q_n \subseteq Y$ for $n \in \omega$ so that:
\begin{itemizz}
\item[1.] $P_{(\,)} = X$ and $Q_0 = Y$.
\item[2.] $P_s$ is closed in $X$ and $Q_n$ is closed in $Y$.
\item[3.] $Q_n = \bigcap\{f(P_s) : \lh(s) = n\}$.
\item[4.] $P_{s\cat 0}$ and $P_{s\cat 1}$ are disjoint subsets of $P_s$.
\item[5.] $\nu(f(P_s) \setminus f(P_{s\cat i})) \le 6^{-n-1}$ 
when $\lh(s) = n$ and $i = 0,1$.
\item[6.] $Q_{n+1} \subseteq Q_n$ and
$\nu(Q_n \backslash Q_{n+1}) \le 2^{n+1} \cdot 6^{-n-1} =  3^{-n-1}$.
\item[7.] $\EEE_\mu(P_s | f)(y) > 0$ for $\nu$--a.e.\@ $y \in f(P_s)$.
\end{itemizz}
Assuming that this can be done, let $Q = \bigcap_n Q_n$. 
$Q \subseteq f(P_s)$ for all $s \in 2^{<\omega}$,
so for $t \in 2^\omega$, let $P_t = f\iv(Q) \cap \bigcap_n P_{t\res n}$.
Then the $P_t$ are disjoint and $f(P_t) = Q$ for all $t$.
Also, $\mu(Q) \ge 1 - 1/3 - 1/9 - 1/27 - \cdots  = 1/2$.
Let $\dom(\varphi) = \bigcup_t P_t$, with $\varphi(x) = t$ for
$x \in P_t$.

Now, to do the construction, note first that (6) follows from (3)(4)(5).
We proceed by induction on $\lh(s)$, using (7) to accomplish the splitting.
For $\lh(s) = 0$, (1)(2)(3)(7) are trivial, since 
$\EEE(X | f)(y) = 1$ for a.e.\@ $y \in Y$.
Now fix $s$ with $\lh(s)  = n$.  We obtain $P_{s\cat 0}$ and $P_{s\cat 1}$ 
by applying Lemma \ref{lemma-nowhere-slim-split},
with the $X,Y$ there replaced by $P_s, f(P_s)$; but then we must
replace $\nu$ by $\lambda := (\mu\res P_s)\, (f\res P_s)\iv$ on $f(P_s)$.
Let $\varphi = \EEE_\mu(P_s | f)$;
then, by (7) for $P_s$, $\varphi(y) > 0$ for $\nu$--a.e.\@ $y \in f(P_s)$; also
$\varphi(y) = 0$ for a.e.\@ $y \notin f(P_s)$, and
$\int_A \varphi(y) \, d\nu(y) = \mu( f\iv(A) \cap P_s) = \lambda(A)$
for all measurable $A \subseteq f(P_s)$.
Fix $\delta > 0$ such that
$\nu(\{y \in f(P_s) :
\varphi(y) <\nobreak\delta \})\allowbreak\le  6^{-n-1}/2 $.
Now apply Lemma \ref{lemma-nowhere-slim-split} to get closed
$P_{s\cat 0}, P_{s\cat 1}$ satisfying (4) with
$K_s := f(P_{s\cat 0}) =  f(P_{s\cat 1})$ so that,
for $i = 0,1$,
$\;  \EEE_{\mu \res P_s}(P_{s\cat i} \mid f \res P_s)(y) > 0$ for
$\lambda$--a.e.\@ $y \in K_s$,
and $\lambda(f(P_s) \backslash K_s) < \delta \cdot 6^{-n-1}/2$.
Now, by Lemma \ref{lemma-realate-exp},
$\EEE_\mu(P_{s\cat i}|f) =
\varphi \cdot \EEE_{\mu \res P_s} (P_{s\cat i} \mid f\res P_s)$,
which yields (7) for $P_{s\cat i}$.
To obtain (5), let $A = f(P_s) \setminus K_s$.
we need $\nu(A) \le 6^{-n-1}$,
and we have
$ \int_A \varphi(y) \, d\nu(y) = \lambda(A) < \delta \cdot 6^{-n-1}/2$.
Let $A = A' \cup A''$, where $\varphi < \delta$ on $A'$ and 
$\varphi \ge \delta$ on $A''$.
Then $\nu(A') \le 6^{-n-1}/2$ and
$\nu(A'') \le (1/\delta) \int_{A''} \varphi(y) \, d\nu(y) \le 6^{-n-1}/2$,
so $\nu(A) \le 6^{-n-1}$.
\end{proofof}

\begin{corollary}
\label{cor-tight-to-separable}
Suppose that $X,Y$ are compact,
$f: X \onto Y$ is weakly $\cccc$--tight,
and $\mu$ is a Radon measure on $X$,
with $\nu = \mu f\iv$ atomless and separable.
Then $\mu$ is separable.
\end{corollary}
\begin{proof}
$X$ is simple with respect to $f,\mu$, by Theorem \ref{thm-tight-to-simple},
which implies that $\ma(\mu)$ is a countable disjoint sum
of separable measure algebras.
\end{proof}

\begin{proofof}{Theorem \ref{thm-dis-ms}}
Assume that $\mu$ is a non-separable Radon measure on $X$; we shall
derive a contradiction.  By subtracting the point masses,
we may assume that $\mu$ is atomless.

First, fix a compact metric $Z$ and a $g : X \onto Z$ such
that $\mu g\iv$ is atomless.
This is easily done by an elementary submodel argument.
More concretely, one can assume that $X \subseteq [0,1]^\kappa$;
then $g = \pi^\kappa_d$ for a suitably chosen countable $d \subseteq \kappa$.
We construct $d$ as $\bigcup_i d_i$, where the $d_i$ are finite
and non-empty and $d_0 \subseteq d_1 \subseteq \cdots$.
Given $d_i$, we have the space $Z_i  = \pi^\kappa_{d_i}(X)$,
with measure $\nu_i = \mu (\pi^\kappa_{d_i})\iv$.
Let $\{F_i^\ell : \ell \in \omega\}$ be a family of closed non-null subsets
of $Z_i$ which is dense in the measure algebra, and make sure
that for each $\ell$, there is some $j > i$ such that
$Z_j$ contains a closed set $K \subseteq (\pi^{d_j}_{d_i})\iv(F_i^\ell)$
with $\nu_j(K) / \mu_i(F_i^\ell) \in (1/3, 2/3)$.

Let $f : X \onto Y$ be weakly $\cccc$--tight, where $Y$ is metric and
$f$ is finer than $g$. 
We then have $\Gamma \in C(Y,Z)$ such that $g = \Gamma \circ f$,
so $\mu g\iv = (\mu f\iv) \Gamma\iv$, so $\mu f\iv$ is atomless.
Also, $\mu f\iv$ is separable because $Y$ is metric,
contradicting Corollary \ref{cor-tight-to-separable}.
\end{proofof}

\section{Inverse Limits}
\label{sec-inv}
Some compacta built as inverse limits in $\omega_1$ steps
are dissipated.
We avoid explicit use of the inverse limit by
viewing $X$ as a subspace of some $M^{\omega_1}$, so
the bonding maps in the inverse limit will be the projection maps.

\begin{definition}
For any space $M$ and ordinals $\alpha \le \beta$:
$\pi^\beta_\alpha: M^\beta \onto M^\alpha$ denotes the natural
projection.
\end{definition}

\begin{theorem}
\label{thm-inv-lim}
Let $M$ be compact metric, and suppose that $X$ is a closed subset
of $M^{\omega_1}$.  Let $X_\alpha = \pi^{\omega_1}_\alpha(X)$.
Assume that for each $\alpha< \omega_1$, the map
$\pi^{\alpha+1}_\alpha \res X_{\alpha+1}: X_{\alpha+1} \onto X_\alpha$
is tight.  Then
\begin{itemizz}
\item[1.] For each $\alpha < \beta \le \omega_1$, the map
$\pi^{\beta}_\alpha \res X_{\beta}: X_{\beta} \onto X_\alpha$
is tight.  
\item[2.] $X$ is dissipated.
\end{itemizz}
\end{theorem}
\begin{proof}
For (1), fix $\alpha$ and induct on $\beta$.
For successor stages, use Lemma \ref{lemma-tight-comp}.
For limit $\beta > \alpha$, use the fact that if
$P_0, P_1$ are disjoint closed subsets of $X_\beta$,
then there is a $\delta$ with $\alpha < \delta < \beta$ and
$\pi^\beta_\delta(P_0) \cap \pi^\beta_\delta(P_1) = \emptyset $.

For (2), observe that given
$g : X \onto Z$, with $Z$ metric, there is an $\alpha < \omega_1$
with $\pi^{\omega_1}_\alpha \res X$ finer than $g$.
Now, use the fact that all $\pi^{\omega_1}_\beta \res X$ are tight.
\end{proof}

The proof of (2) did not actually require all
$\pi^{\omega_1}_\beta \res X$ to be tight;
we only needed unboundedly many.
More generally, the definition of ``dissipated'' requires the
family of tight maps
to be unbounded, but it does not necessarily contain a club, although
it does contain a club in the ``natural'' examples of dissipated spaces.
We first point out
an example where the tight maps do not contain a club.
Then we shall formulate precisely what ``contains a club'' means.

\begin{example}
\label{ex-noclub}
There is a closed $X \subseteq 2^{\omega_1}$ such that,
setting  $X_\alpha = \pi^{\omega_1}_\alpha(X)$:
\begin{itemizz}
\item[a.] $X$ is dissipated
\item[b.] For all $\alpha < \omega_1$, 
$\pi^{\omega_1}_\alpha\res X : X \onto X_\alpha$ is tight iff
$\alpha$ is not a limit ordinal.
\end{itemizz}
\end{example}
\begin{proof}
First note that $(b) \to (a)$ because 
whenever $g : X \to Z$, with $Z$ metric,
there is always an $\alpha < \omega_1$ with
$\pi^{\omega_1}_\alpha\res X \le g$.  Then
$\pi^{\omega_1}_{\alpha+1}\res X  \le \pi^{\omega_1}_\alpha\res X \le g$ and
$\pi^{\omega_1}_{\alpha+1} \res X$ is tight.

To prove $(b)$, we use a standard inverse limit construction,
building $X_\alpha$ by induction on $\alpha$.
We shall have:
\begin{itemizz}
\item[1.] $X_\alpha$ is a closed subset of $2^\alpha$ for all
$\alpha \le \omega_1$, and $X = X_{\omega_1}$.
\item[2.] $X_\alpha = \pi^\beta_\alpha(X_\beta)$
whenever $\alpha \le \beta \le \omega_1$.
\item[3.] $X_\alpha = 2^\alpha$ for $\alpha \le \omega$.
\item[4.] For $\alpha < \omega_1$:
$X_{\alpha + 1} = X_\alpha \times \{0\} \cup F_\alpha \times \{1\}$,
where $F_\alpha$ is a closed subset of $X_\alpha$.
\item[5.] $F_\gamma$ is a perfect set for all limit $\gamma < \omega_1$.
\item[6.] $\pi^\alpha_\delta(F_\alpha)$ is finite
whenever $\delta < \alpha < \omega_1$.
\item[7.] Whenever $\delta < \alpha < \omega_1$ and $\delta$ is 
a successor ordinal, there is an $n$ with $0 < n < \omega$
such that
$\pi^{\alpha + n}_{\delta + 1}(F_{\alpha + n}) = F_\delta \times \{0,1\}$.
\end{itemizz}
Conditions (1)(2) imply that $X_\gamma$, for limit $\gamma$,
is determined by the $X_\alpha$ for $\alpha < \gamma$;
then, by (4), the whole construction is determined by the
choice of the $F_\alpha \subseteq  X_\alpha$;
as usual, in stating (4), we are identifying
$2^{\alpha + 1}$ with $2^\alpha \times \{0,1\}$.
By (3),  $F_\alpha = X_\alpha$ when $\alpha < \omega$.
By (6), $F_\alpha$ is finite for successor $\alpha$.
Conditions (1)--(6) are sufficient to verify $(b)$ of the theorem,
but (7) was added to ensure that the construction can be carried out.
Using (7), it is easy to construct $F_\gamma$ for limit $\gamma$
to satisfy (5)(6), and (7) itself is easy to ensure by a standard
enumeration argument, since there are no further restrictions
on the finite sets $F_{\alpha + n} \subseteq X_{\alpha + n}$ when $n > 0$.

To verify $(b)$:  If $\alpha < \omega_1$ is a limit ordinal,
then (4)(5) guarantee that 
$\pi^{\omega_1}_\alpha\res X : X \onto X_\alpha$ is not tight.
Now, fix a successor $\alpha < \omega$. 
We prove by induction that 
$\pi^{\beta}_\alpha\res X_\beta : X_\beta \onto X_\alpha$ is tight
whenever $\alpha \le \beta \le \omega_1$.
This is trivial when $\beta = \alpha$.  If 
$\beta > \alpha$ is a limit ordinal and $\pi^{\beta}_\alpha\res X_\beta $
fails to be tight, then we have disjoint closed $P_0,P_1 \subset X_\beta$
with $Q = \pi^{\beta}_\alpha(P_0) = \pi^{\beta}_\alpha(P_1)$ and $Q$ not
scattered; but then there is a $\delta$ with
$\beta > \delta > \alpha$ such that
$\pi^{\beta}_\delta(P_0) \cap \pi^{\beta}_\delta(P_1) = \emptyset$,
and then the $\pi^{\beta}_\delta(P_i)$ refute the tightness of
$\pi^\delta_\alpha$.

Finally, assume that $\alpha \le \beta < \omega_1$ and that
$\pi^{\beta}_\alpha\res X_\beta $ is tight.  We shall prove that
$\pi^{\beta + 1}_\alpha\res X_{\beta + 1} $ is tight. 
If $\beta$ is a successor, we note that
$\pi^{\beta + 1}_\beta\res X_{\beta + 1} $ is tight
because $F_\beta$ is finite, so that
$\pi^{\beta + 1}_\alpha\res X_{\beta + 1} =
\pi^{\beta}_\alpha\res X_\beta \circ
\pi^{\beta + 1}_\beta\res X_{\beta + 1} 
$ is tight by
Lemma \ref{lemma-tight-comp}.
Now, assume that $\beta$ is a limit (so $\alpha < \beta$) and that 
$\pi^{\beta + 1}_\alpha\res X_{\beta + 1} $ is not tight. 
Fix disjoint closed $P_0,P_1 \subset X_{\beta+1}$
with $Q = \pi^{\beta+1}_\alpha(P_0) = \pi^{\beta+1}_\alpha(P_1)$ and $Q$ not
scattered.
Since $\pi^\beta_\alpha(F_\beta)$ is finite, we may shrink $Q$ and
the $P_i$ and assume that $Q \cap \pi^\beta_\alpha(F_\beta) = \emptyset$.
Then $\pi^{\beta+1}_\beta(P_i) \cap F_\beta = \emptyset$, so that
$\pi^{\beta+1}_\beta(P_0) \cap \pi^{\beta+1}_\beta(P_1) = \emptyset$,
and the $\pi^{\beta+1}_\beta(P_i)$ contradict the tightness 
of $\pi^{\beta}_\alpha\res X_\beta $.
\end{proof}

There are various equivalent ways to formulate ``contains a club'';
the following is probably the quickest to state:

\begin{definition}
The compact $X$ is \emph{wasted} iff whenever $\theta$ is a suitably
large regular cardinal and $M \prec H(\theta)$ is countable and contains
$X$ and its topology, the natural
evaluation map $\pi_M : X \to  [0,1]^{C(X, [0,1]) \cap M} $ is tight.
\end{definition}

For the $X$ of Example \ref{ex-noclub}, no $\pi_M$ is tight,
since $\pi_M$ is equivalent to $\pi^{\omega_1}_\gamma$, where
$\gamma = \omega_1 \cap M$.
The $X$ of Theorem \ref{thm-inv-lim} is wasted, as is
every compact LOTS.
A notion intermediate between ``dissipated'' and ``wasted''
is obtained by requiring $\pi_M $ to be tight for
a stationary set of $M \prec H(\theta)$.

In Theorem \ref{thm-inv-lim}: 
since $X_{\alpha+1}$ and   $X_{\alpha}$ are compact
metric, the assumption that $\pi^{\alpha+1}_\alpha $ is tight is equivalent
to saying that
$\{y \in X_\alpha : 
| (\pi^{\alpha + 1}_\alpha)\iv\{y\} \cap X_{\alpha + 1} | > 1 \}$
is countable (see Theorem \ref{thm-tight-metric}).
In the constructions of \cite{Fed1, HK2, HK3},
this set is actually a singleton.
In some cases, the spaces are also
\emph{minimally generated} in the sense 
Koppelberg \cite{Kop} and Dow \cite{Dow}:

\begin{definition}
Let $X,Y$ be compact.  Then $f : X \onto Y$ is \emph{minimal}
iff $|f\iv\{y\}| = 1$ for all $y \in Y$ except for one $y_0$,
for which $|f\iv\{y_0\}| = 2$.
\end{definition}

We remark that this is the same as minimality in the sense that
if $f = g \circ h$, where $h : X \onto Z$ and $g : Z \onto Y$,
then either $g$ or $h$ is a bijection.
Clearly, every minimal map is tight.

\begin{definition}
$X$ is \emph{minimally generated} iff $X$ is a closed subspace of
some $2^\rho$, where, setting
$X_\alpha = \pi^{\rho}_\alpha(X)$, all the maps
$\pi^{\alpha+1}_\alpha \res X_{\alpha+1}: X_{\alpha+1} \onto X_\alpha$,
for $\alpha < \rho$, are minimal.
\end{definition}

Examples of such spaces are the 
Fedor\v cuk S-space \cite{Fed1}, obtained under $\diamondsuit$
(here, $\rho = \omega_1$), and the Efimov spaces
obtained by Fedor\v cuk \cite{Fed2} and Dow \cite{Dow},
where  $\rho > \omega_1$.

Clearly, if $\rho = \omega_1$, then $X$ must be dissipated
by Theorem \ref{thm-inv-lim}, but this need not be true for $\rho > \omega_1$.
For example, if $A(\aleph_1)$ is the 1-point compactification of a
discrete space of size $\aleph_1$, and $X = A(\aleph_1) \times 2^\omega$,
then $X$ is not $\aleph_1$--dissipated by Lemma \ref{lemma-prod-bad},
but $X$ is minimally generated, with $\rho = \omega_1 + \omega$.

Note that if we weaken ``tight'' to ``3--tight'' in Theorem \ref{thm-inv-lim},
we get nothing of any interest in general. 
In fact, if $M = 2 = \{0,1\}$ and each
$X_\alpha = M^\alpha$, then
all $\pi^{\alpha+1}_\alpha \res X_{\alpha+1}$ are 3--tight,
but $X$ is not weakly $\cccc$--dissipated by Theorem 
\ref{thm-bad-product}.
However, one can in some cases use an inverse limit construction
build a space which is $\aleph_0$--dissipated:

\begin{proofof}{Proposition \ref{prop-L}}
We modify the standard construction of a compact L-space under CH,
following specifically the details in \cite{Kun2};
similar constructions are in  Haydon \cite{Haydon} and Talagrand \cite{Tal}.
So, $X$ will be a closed subset of $2^{\omega_1}$.

We inductively define $X_\alpha \subseteq 2^\alpha$,
for $\omega \le \alpha \le \omega_1$,
along with an atomless Radon probability measure $\mu_\alpha$ on $X_\alpha$
such that the support of $\mu_\alpha$ is all of $X_\alpha$.
Let $X_\omega = 2^\omega$ with $\mu_\omega$
the usual product measure.
The measures will all cohere, in the sense that 
$\mu_\alpha = \mu_\beta\, (\pi^\beta_\alpha)\iv$ whenever
$\alpha < \beta$.
Along with the measures, we choose a countable family $\FF_\alpha$
of closed $\mu_\alpha$--null subsets of $X_\alpha$
and a specific closed nowhere dense non-null $K_\alpha \subseteq X_\alpha$.
When $\alpha < \beta < \omega_1$, $\FF_\beta$
will contain $(\pi^\beta_\alpha)\iv(F)$ for all $F \in \FF_\beta$,
along with some additional sets.
Since $\FF_\alpha$ is countable, we can choose a perfect 
$C_\alpha \subseteq K_\alpha$ such that 
$\mu_\alpha(C_\alpha) > 0$, $C_\alpha$ is the support of
$\mu_\alpha\res C_\alpha$, and
$C_\alpha \cap F = \emptyset$ for all $F \in \FF_\alpha$.
Then we let
$X_{\alpha + 1} = X_\alpha \times \{0\} \cup C_\alpha \times \{1\}$.
In the construction of \cite{Kun2}, 
$\mu_{\alpha + 1}$ can be chosen arbitrarily to satisfy
$\mu_\alpha = \mu_{\alpha + 1}\, (\pi^{\alpha + 1}_\alpha)\iv$,
as long as all non-empty open subsets of $ C_\alpha \times \{1\}$
have positive measure;
there is some flexibility here in distributing the measure on
$C_\alpha$ among its copies 
$C_\alpha \times \{0\}$ and  $ C_\alpha \times \{1\}$.
In particular,  depending on the choices made,
the final measure $\mu = \mu_{\omega_1}$ on $X = X_{\omega_1}$
may be separable or non-separable.
In any case, \cite{Kun2} shows that, assuming CH, one may choose the
$\FF_\alpha$ and $K_\alpha$ appropriately to guarantee
$X$ is an L-space and that the ideals
of null subsets, meager subsets, and separable subsets
all coincide.

Now, always choose $\mu_{\alpha+1}$ such that 
$\mu_{\alpha+1}( C_\alpha \times \{0\} ) = 0$.  This will guarantee
that $\mu$ on $X$ is separable, with $\ma(\mu)$ isomorphic to
$\ma(\mu_\omega)$ via $(\pi^{\omega_1}_\omega)^*$.
Also, put the set $C_\alpha \times \{0\}$ into $\FF_{\alpha+1}$.
Then, for all $x \in X_\omega$,
$(\pi^{\omega_1}_\omega)\iv\{x\}$ is scattered (as is easy to verify),
and hence countable (since $X$ is HL).
But then  $\pi^{\omega_1}_\omega \res X: X \onto X_\omega$
is $\aleph_1$--tight, so that $X$ is
$\aleph_1$--dissipated by Lemma \ref{lemma-one-enough}.
\end{proofof}

We remark that by Theorem \ref{thm-dis-ms}, we know that
the $\mu$ of Proposition \ref{prop-L}
must be separable, so it was natural to make
$\ma(\mu)$ isomorphic to $\ma(\mu_\omega)$ in the construction.

\section{Absoluteness}
\label{sec-ab}
We shall prove here that tightness is absolute.
This can then be applied in forcing arguments, but the absoluteness itself has
nothing at all to do with forcing; it is just a fact about
transitive models of ZFC, and is related to the absoluteness
of $\Pi^1_1$ statements.  Since we never need absoluteness of 
$\Pi^1_2$ (Shoenfield's Theorem), we do not need the
models to contain all the ordinals.
So, we consider arbitrary transitive models $M,N$ of ZFC
with $M \subseteq N$.
If in $M$, we have compacta $X,Y$ and $f: X \to Y$, we want to show
that $f$ is tight in $M$ iff $f$ is tight in $N$.

To make this discussion precise, we must, in $N$,
replace $X,Y$ by the corresponding compact spaces
$\widetilde X,\widetilde Y$.
This concept was  
described by Bandlow \cite{BAND} (and later in \cite{Dow,DF,DK2,HK3}),
and is defined as follows:

\begin{definition}
\label{def-tilde}
Let $M \subseteq N$ be transitive models of ZFC.
In $M$, assume that $X$ is compact.
Then $\widetilde X$ denotes the compactum in $N$ characterized by:
\begin{itemizz}
\item[1.] $X$ is dense in $\widetilde X$.
\item[2.] Every $\varphi \in C(X, [0,1]) \cap M$
extends to a $\widetilde \varphi \in C(\widetilde X, [0,1])$ in  $N$.
\item[3.] The functions $\widetilde \varphi$ \textup(for $\varphi \in M$\textup)
separate the points of $\widetilde X$.
\end{itemizz}
If, in $M$,  $X,Y$ are compact and $f \in C(X,Y)$, then 
in $N$, $\widetilde f \in C(\widetilde X,\widetilde Y)$
denotes the \textup(unique\textup) continuous extension of $f$.

In forcing, $\nameX$ denotes the $\widetilde X$ of $V[G]$,
while $\check X$ denotes the $X$ of $V[G]$.
\end{definition}

\begin{theorem}
\label{thm-absol}
Let $M \subseteq N$ be transitive models of ZFC.
In $M$, assume that $X,Y$ are compact, $K$ is compact metric,
and $f: X \to Y$.  Then the following are equivalent:
\begin{itemizz}
\item[1.] In $M$:  There is a $K$--loose function for $f$.
\item[2.] In $N$:  There is a $\widetilde K$--loose function for $\widetilde f$.
\end{itemizz}
\end{theorem}
\begin{proof}
For $(1) \to (2)$, just observe that if in $M$, we have
$\varphi, Q$ satisfying Definition \ref{def-K-tight} (of $K$--loose),
then $\widetilde \varphi, \widetilde Q$ satisfy 
Definition \ref{def-K-tight} in $N$.

For $\neg(1) \to \neg(2)$, we shall define a partial order $\TTT$ in $M$.
We shall then prove that  $\neg(1)$ implies the well-founded of $\TTT$ in $M$,
while the well-founded of $\TTT$ in $N$ implies $\neg(2)$.
The result then follows by the absoluteness of well-foundedness.

As in the proof of Theorem \ref{thm-metric-cover},
let $H = [0,1]^\omega$, and assume that $K \subseteq H$.
Then the existence of a $K$--loose function is equivalent
to the existence of a $\varphi \in C(X,H)$ such that
for some non-scattered $Q \subseteq Y$ we have $\psi(f\iv\{y\}) \supseteq K$
for all $y \in Q$.

$\TTT$ is a tree of finite sequences, ordered by extension.
$\TTT$ contains the empty sequence and all non-empty
sequences 
\[
\big \langle (\EE_0, \psi_0),  (\EE_1, \psi_1),  \ldots , (\EE_{n-1},
\psi_{n-1})   \big \rangle
\]
satisfying:
\begin{itemizz}
\item[a.] Each $\psi_i \in C(X,H)$.
\item[b.] Each $\EE_i$ is a disjoint family of $2^i$ non-empty
closed subsets of $Y$.
\item[c.]  Whenever $y \in E \in \EE_i$ and $z \in K$:
$d(z, \psi_i(f\iv\{y\})) \le 2^{-i}$.
\item[d.] When $i+1 < n$:  $d(\psi_i, \psi_{i+1}) \le 2^{i-1}$,
and each $E \in \EE_i$ has exactly two subsets in $\EE_{i+1}$.
\end{itemizz}
In $M$, if $\TTT$ is not well-founded and 
$ \langle (\EE_0, \psi_0),  (\EE_1, \psi_1), \ldots \rangle$
is an infinite path through $\TTT$, then we get 
$\varphi = \lim_i \psi_i \in C(X,H)$ using (a)(d)
and $Q = \bigcap_i \bigcup \EE_i$, which is a non-scattered subset
of $Y$ using (b)(d), and (c)(d) implies that 
$\varphi( f\iv \{y\} ) \supseteq K$ for all $y \in Q$, so (1) holds.

Now, suppose, in $N$, that we have $Q, \varphi$ for which (2) holds;
then we construct a path through $\TTT$.
To obtain the $\psi_i$ (all of which must be in $M$),
use the fact that $\{\widetilde \psi : \psi \in C(X,H)^M\}$
is dense in $C(\widetilde X, \widetilde H)$.
Likewise each $E \in \EE_i$ will be a closed set in $M$
such that $\widetilde E \cap Q$ is not scattered.
\end{proof}

Note that Theorem \ref{thm-absol} says that the existence of
the $\varphi$ and $Q$ described in the proof Theorem \ref{thm-metric-cover}
is absolute.  The corresponding ``absoluteness version'' of
Theorem \ref{thm-metric-embed} is false.
For example, suppose that in $V$, we have $X = Y \times K$,
where $X,Y,K$ are compact and non-scattered, and in addition,
$K$ has no non-trivial convergent $\omega$--sequences.
Then clearly in $V$, there can be no
perfect $Q \subseteq Y$ and 1-1 map $i : Q \times (\omega+1) \to X$
such that $f(i(q,u)) = q$ for all $(q,y) \in Q\times (\omega+1)$,
whereas if $V[G]$ collapses enough cardinals, it will
contain such a $Q,i$.

An application of the absoluteness result in Theorem
\ref{thm-absol} is:

\begin{proofof}{Theorem \ref{thm-omega-loose}}
Assume that in the universe, $V$: $X$ and $Y$ are compact, $f: X \to Y$,
and we have an infinite loose family $\{P_i : i \in \omega\}$.
Let $V[G]$ be any forcing extension of $V$ which makes the
weights of $X$ and $Y$ countable, so that in $V[G]$,
we still have
$f: \widetilde X \to \widetilde Y$ and a loose family
$\{\widetilde{P_i}: i \in \omega\}$, but 
$\widetilde X$ and $\widetilde Y$ are now compact metric, so that
Theorem \ref{thm-metric-cover} gives us an $(\omega + 1)$--loose function
in $V[G]$.  Hence, by absoluteness, there is one in $V$.
\end{proofof}

A direct proof of this can be given without forcing, but it seems
quite a bit more complicated, since one must embed into the proof
the method of Suslin used in proving Lemma \ref{lemma-suslin};
one cannot just quote Suslin's theorem, since the spaces are not Polish.
Theorem \ref{thm-omega-loose} is needed for the $\kappa = \omega$
part of:

\begin{corollary}
\label{cor-absol}
Fix $\kappa \le \omega$.
Let $M,N$ be transitive models of ZFC, with $M \subseteq N$.
Assume that in $M$ we have $X,Y,f$ with
$X,Y$ compact and $f: X \to Y$. 
Then $M \models \mbox{``$f: X \to Y$ is $\kappa$--tight''}$  iff
$N \models \mbox{``$\widetilde f: \widetilde X \to \widetilde Y$
is $\kappa$--tight''}$.
\end{corollary}

Of course, the $\leftarrow$ direction is trivial, and holds for 
all $\kappa$ if we rephrase Definition \ref{def-tight} appropriately
so that $\kappa$ is not required to be a cardinal (since
``cardinal'' is not absolute).
That is, if in $M$, we have a loose family
$\{P_\alpha : \alpha < \kappa\}$, then 
$\{\widetilde{P_\alpha} : \alpha < \kappa\}$ is loose in $N$.
For a version of Corollary \ref{cor-absol} for $\kappa = \cccc$,
we use the notion of ``weakly $\cccc$--tight'' from
Definition \ref{def-weak-tight}.

\begin{corollary}
\label{cor-absol-c}
Fix $\kappa \le \omega$.
Let $M,N$ be transitive models of ZFC, with $M \subseteq N$.
Assume that in $M$ we have $X,Y,f$ with
$X,Y$ compact and $f: X \to Y$. 
Then $M \models \mbox{``$f: X \to Y$ is weakly $\cccc$--tight''}$  iff
$N \models \mbox{``$\widetilde f: \widetilde X \to \widetilde Y$
is weakly $\cccc$--tight''}$.
\end{corollary}

\end{document}